\documentclass{cmslatex}
\usepackage[english]{babel}
\usepackage[utf8]{inputenc}
\usepackage[squaren]{SIunits}
\usepackage{latexsym}
\usepackage{amssymb}
\usepackage{enumerate}
\usepackage{amsmath}
\usepackage{pgfplots}
\usepackage[margin=1em]{subcaption}
\usepackage{tikz}
\usepackage{graphicx}
\usepackage{url}
\usepackage{booktabs}
\usepackage{hyperref}
\usepackage{cleveref}

\renewcommand{\theequation}{\arabic{section}.\arabic{equation}}

\crefname{equation}{Eq.}{Eqs.}
\creflabelformat{equation}{#2(#1)#3}
\crefname{section}{Section}{Sections}
\crefname{figure}{Figure}{Figures}
\crefname{table}{Table}{Tables}

\graphicspath{{./}{figures/}}
\usetikzlibrary{decorations.pathreplacing}
\usetikzlibrary{calc}
\usetikzlibrary{arrows}
\usetikzlibrary{positioning}
\usetikzlibrary{shapes.geometric}
\pgfplotsset{
  compat=1.5,
  xmin=0.01,
  xmax=1.00,
}

\makeatletter
\newcommand*{\dif}{\@ifnextchar^{\DIfF}{\DIfF^{}}}
\def\DIfF^#1{\mathop{\mathrm{\mathstrut d}}\nolimits^{#1}\gobblesp@ce}
\def\gobblesp@ce{\futurelet\diffarg\opsp@ce}
\def\opsp@ce{%
  \let\DiffSpace\!%
  \ifx\diffarg(%
    \let\DiffSpace\relax
  \else
    \ifx\diffarg[%
      \let\DiffSpace\relax
    \else
      \ifx\diffarg\{%
        \let\DiffSpace\relax
      \fi\fi\fi\DiffSpace}
\makeatother

\newcommand*{\vct}[1]{\ensuremath{\boldsymbol{#1}}}
\newcommand*{\bigo}[1]{\ensuremath{\mathcal{O}\left(#1\right)}}
\newcommand*{\txin}{\ensuremath{\textup{in}\ }}
\newcommand*{\txon}{\ensuremath{\textup{on}\ }}
\newcommand*{\smallo}[1]{\ensuremath{o\left(#1\right)}}
\newcommand*{\eint}[1]{\ensuremath{\int_{-\infty}^\infty #1\dif z}}

\newcommand*{\zlimpm}{\ensuremath{\lim_{z\to\pm\infty}}}

\renewcommand*{\div}{\boldsymbol\nabla\cdot}
\newcommand*{\divs}{\boldsymbol\nabla_{\text s}\cdot}
\renewcommand*{\grad}{\boldsymbol\nabla}
\newcommand*{\grads}{\boldsymbol\nabla_{\text s}}
\newcommand*{\lapls}{\boldsymbol\Delta_{\text s}}
\newcommand*{\lapl}{\boldsymbol\Delta}

\newcommand*{\ndot}{\vct n\cdot}
\newcommand*{\dotn}{\cdot\vct n}
\newcommand*{\e}[1]{\ensuremath{\cdot 10^{-#1}}}

\begin{document}

\title{Analysis of the diffuse-domain method for solving PDEs in
  complex geometries\thanks{Submitted to Communications in Mathematical
    Sciences}}

\author{Karl Yngve Lervåg\thanks{Department of Energy and Process Engineering,
    Norwegian University of Science and Technology, NO-7491 Trondheim,
    Norway.  SINTEF Energy Research, P.O.\ Box 4761 Sluppen, NO-7465
    Trondheim, Norway (karl.lervag@sintef.no)}
  \and John Lowengrub\thanks{Department of Mathematics, University of
    California, Irvine, Irvine CA-92697, USA (lowengrb@math.uci.edu).}}

\pagestyle{myheadings}%
  \markboth{\textsc{Analysis of the diffuse-domain method}}%
  {\textsc{Karl Yngve Lervåg and John Lowengrub}}
\maketitle

\begin{abstract}
  In recent work, Li et al.\ (Comm.\ Math.\ Sci., 7:81-107, 2009) developed
  a diffuse-domain method (DDM) for solving partial differential equations in
  complex, dynamic geometries with Dirichlet, Neumann, and Robin boundary
  conditions.  The diffuse-domain method uses an implicit representation of the
  geometry where the sharp boundary is replaced by a diffuse layer with
  thickness $\epsilon$ that is typically proportional to the minimum grid size.
  The original equations are reformulated on a larger regular domain and the
  boundary conditions are incorporated via singular source terms.  The
  resulting equations can be solved with standard finite difference and finite
  element software packages.  Here, we present a matched asymptotic analysis of
  general diffuse-domain methods for Neumann and Robin boundary conditions.
  Our analysis shows that for certain choices of the boundary condition
  approximations, the DDM is second-order accurate in $\epsilon$.  However, for
  other choices the DDM is only first-order accurate.  This helps to explain
  why the choice of boundary-condition approximation is important for rapid
  global convergence and high accuracy.  Our analysis also suggests correction
  terms that may be added to yield more accurate diffuse-domain methods.
  Simple modifications of first-order boundary condition approximations are
  proposed to achieve asymptotically second-order accurate schemes.  Our
  analytic results are confirmed numerically in the $L^2$ and $L^\infty$ norms
  for selected test problems.
\end{abstract}


\section{Introduction}

There are many problems in computational physics that involve solving partial
differential equations (PDEs) in complex geometries.  Examples include fluid
flows in complicated systems, vein networks in plant leaves, and tumours in
human bodies.  Standard solution methods for PDEs in complex domains typically
involve triangulation and unstructured grids.  This rules out coarse-scale
discretizations and thus efficient geometric multi-level solutions.  Also, mesh
generation for three-dimensional complex geometries remains a challenge, in
particular if we allow the geometry to evolve with time.

In the past several years, there has been much effort put into the development
of numerical methods for solving partial differential equations in complex
domains.  However, most of these methods typically require tools
not frequently available in standard finite element and finite difference
software packages.  Examples of such approaches include the extended and
composite finite element methods (e.g.,
\cite{GR07,dolbow09,fries10,duddu11,he11,PRE11,byfut12,bernauer12}), immersed
interface methods (e.g., \cite{LL94,LiIto_2006,SethianShan_2008,li12,wan12}),
virtual node methods with embedded boundary conditions (e.g.,
  \cite{Bedrossian10,Zhu12,Hellrung12}), matched interface and boundary methods
(e.g., \cite{zhou06,zhao09,zhao10,xia11,zhou12}), modified finite
volume/embedded boundary/cut-cell methods/ghost-fluid methods (e.g.,
  \cite{GlimmMarchesinMcBryan_1981,JohansenColella_1998,FedkiwAslamMerrimanOsher_1999,GibouFedkiwChengKang_2002,GibouFedkiw_2005,JiLienYee_2006,MacklinLowengrub_2006,zhong07,MacklinLowengrub_2008,Colellaetal_2008,lui09,Uzgoren09,OSK09,cisternino12,Coco13,Papac10,Papac13,Helgadottir11,Theillard13}).
In another approach, known as the fictitious domain method (e.g.,
  \cite{GlowinskiPanPeriaux_CMAME_1994,GlowinskiPanWellsZhou_JCP_1996,RamiereAngotBelliard_CMAME_2007,lohner07}),
the original system is either augmented with equations for Lagrange multipliers
to enforce the boundary conditions, or the penalty method is used to enforce
the boundary conditions weakly.  See also \cite{Gibou13} for a review of
numerical methods for solving the Poisson equation, the diffusion equation and
the Stefan problem on irregular domains.

An alternate approach for simulating PDEs in complex domains, which does not
require any modification of standard finite element or finite difference
software, is the diffuse-domain method.  In this method, the domain is
represented implicitly by a phase-field function, which is an approximation of
the characteristic function of the domain.  The domain boundary is replaced by
a narrow diffuse interface layer such that the phase-field function rapidly
transitions from one inside the domain to zero in the exterior of the domain.
The boundary of the domain can thus be represented as an isosurface of the
phase-field function.  The PDE is then reformulated on a larger, regular domain
with additional source terms that approximate the boundary conditions.
Although uniform grids can be used, local grid refinement near domain
boundaries improves efficiency and enables the use of smaller interface
thicknesses than are achievable using uniform grids.  A related approach
involves the level-set method \cite{osher88,Sethian99,osher03a} to describe the
implicitly embedded surface and to obtain the appropriate surface operators
(e.g., \cite{greer06}).

The diffuse-domain method (DDM) was introduced by Kockelkoren et al.\
\cite{Kockelkoren03} to study diffusion inside a cell with zero Neumann
boundary conditions at the cell boundary (a similar approach was also used in
\cite{Bueno06b,Bueno06a} using spectral methods).  The DDM was later used to
simulate electrical waves in the heart \cite{Fenton05} and membrane-bound
Turing patterns \cite{Levine05}.  More recently, diffuse-interface methods have
been developed for solving PDEs on stationary \cite{ratz} and evolving
\cite{DD07,dziuk08a,dziuk08b,elliott09a,elliott09,DE12} surfaces.
Diffuse-domain methods for solving PDEs in complex evolving domains with
Dirichlet, Neumann and Robin boundary conditions were developed by Li et al.\
\cite{Li09} and by Teigen et al.\ \cite{Teigen09-b} who modelled bulk-surface
coupling.  The DDM was also used by Aland et al.\ \cite{Aland10} to simulate
incompressible two-phase flows in complex domains in 2D and 3D, and by Teigen
et al.\ \cite{Teigen11} to study two-phase flows with soluble surfactants.

Li et al.\ \cite{Li09} showed that in the DDM there exist several
approximations to the physical boundary conditions that converge asymptotically
to the correct sharp-interface problem.  Li et al.\ presented some numerical
convergence results for a few selected problems and observed that the choice of
boundary condition can significantly affect the accuracy of the DDM.  However,
Li et al.\ did not perform a quantitative comparison between the different
boundary-condition approximations, nor did they estimate convergence rates.
Further, Li et al.\ did not address the source of the different levels of
accuracy they observed for the different boundary-condition approximations.

In the context of Dirichlet boundary conditions, Franz et al.\
\cite{Franz12} recently presented a rigorous error analysis of the DDM for
a reaction-diffusion equation and found that the method converges only with
first-order accuracy in the interface thickness parameter $\epsilon$, which
they confirmed numerically.  Similar results were obtained numerically by
Reuter et al.\ \cite{Reuter12} who reformulated the DDM using an integral
equation solver.  Reuter et al.\ demonstrated that their generalized DDM,
with appropriate choices of approximate surface delta functions, converges
with first-order accuracy to solutions of the Poisson equation with Dirichlet
boundary conditions.

Here, we focus on Neumann and Robin boundary conditions and we present
a matched asymptotic analysis of general diffuse-domain methods in a fixed
complex geometry, focusing on the Poisson equation for Robin boundary
conditions and a steady reaction-diffusion equation for Neumann boundary
conditions.  However, our approach applies to transient problems and more
general equations in the same way as shown in \cite{Li09}.  Our analysis
shows that for certain choices of the boundary condition approximations, the
DDM is second-order accurate in $\epsilon$, which in practice is proportional
to the smallest mesh size.  However, for other choices the DDM is only
first-order accurate.  This helps to explain why the choice of boundary
condition approximation is important for rapid global convergence and high
accuracy.

Further, inspired by the work of Karma and Rappel \cite{Karma98} and Almgren
\cite{Almgren99}, who incorporated second-order corrections in their phase
field models of crystal growth and by the work of Folch et al.\
\cite{Folch1999} who added second-order corrections in phase-field models of
advection, we also suggest correction terms that may be added to yield a more
accurate version of the diffuse-domain method.  Simple modifications of
first-order boundary condition approximations are proposed to achieve
asymptotically second-order accurate schemes.  Our analytic results are
confirmed numerically for selected test problems.

The outline of the paper is as follows.  In \cref{sec:dda} we introduce and
present an analysis of general diffuse-domain methods.  In
\cref{sec:discretization} the numerical methods are described, and in
\cref{sec:results} the test cases are introduced and numerical results are
presented and discussed.  We finally give some concluding remarks in
\cref{sec:conclusion}.

\section{The diffuse-domain method}
\label{sec:dda}

The main idea of the DDM is to extend PDEs that are defined inside complex and
possibly time-dependent domains into larger, regular domains.  As a model
problem, consider the Poisson equation in a domain $D$,
\[
  \lapl u = f,
\]
with Neumann or Robin boundary conditions.  As shown in Li et al.\ \cite{Li09},
the results for the Poisson equation can be used directly to obtain
diffuse-domain methods for more general second-order partial differential
equations in evolving domains.

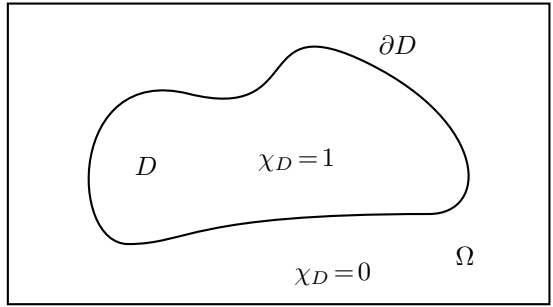
\begin{figure}[tbp]
  \centering
  \begin{tikzpicture}[scale=0.8]
    \draw[thick] (0,0) rectangle(9,5);
    \draw[thick] (2,1)
    .. controls (1,1) and (1,4) .. (3,3.5)
    .. controls (5,3) and (4,5) .. (6,4)
    node[above right] {$\partial D$}
    .. controls (8,3) and (8,1.5) .. (7,1.5)
    .. controls (3,1.5) and (3,1) .. (2,1);
    \node at (2.3,2.3) {$D$};
    \node at (7.6,0.8) {$\Omega$};
    \node at (4.8,2.4) {$\chi_D=1$};
    \node at (5.4,0.5) {$\chi_D=0$};
  \end{tikzpicture}
  \caption{A complex domain $D$ covered by a larger, regular domain $\Omega$.}
  \label{fig:ddadomain}
\end{figure}

The DDM equation is defined in a larger, regular domain $\Omega\supset D$ as
\begin{equation}
  \div(\phi\grad u) + \text{BC} = \phi f,
  \label{ddm}
\end{equation}
see \cref{fig:ddadomain}.  Here $\phi$ approximates the characteristic function
of $D$,
\[
  \chi_D = \begin{cases}
    1 & \text{if $x\in D$,} \\
    0 & \text{if $x\notin D$,}
  \end{cases}
\]
and BC is chosen to approximate the physical boundary condition, cf.\
\cite{Li09}.  This typically involves diffuse-interface approximations of the
surface delta function.  A standard approximation of the characteristic
function is the phase-field function,
\begin{equation}
  \chi_D \simeq \phi(\vct x,t) = \frac{1}{2} \left( 1
  - \tanh \left( \frac{3r(\vct x,t)}{\epsilon} \right) \right).
  \label{eq:characteristic}
\end{equation}
Here $\epsilon$ is the interface thickness and $r(\vct x,t)$ is the
signed-distance function with respect to $\partial D$, which is taken to be
negative inside $D$.

As Li et al.\ \cite{Li09} described, there are a number of different choices
for BC in \cref{ddm}.  For example, in the Neumann case, where $\ndot\nabla
u = g$ on $\partial D$, one may take:
\[
  \text{BC} =
    \begin{cases}
      \text{BC1} = g|\nabla\phi|,
      \quad\text{or}\\
      \text{BC2} = \epsilon g|\nabla\phi|^2.
    \end{cases}
\]
In the Robin case, where $\ndot\nabla u = k(u-g)$ on $\partial D$, one may use
analogous approximations:
\begin{equation}
  \text{BC} =
    \begin{cases}
      \text{BC1} = k(u - g) |\nabla\phi |,
      \quad\text{or} \\
      \text{BC2} = \epsilon k(u - g) |\nabla\phi|^2.
    \end{cases}
  \label{eq:bc_robin}
\end{equation}
Note that the terms $|\nabla\phi|$ and $\epsilon|\nabla\phi|^2$ approximate the
surface delta function.  Following Li et al.\ \cite{Li09} we assume that $g$ is
extended constant in the normal direction off $\partial D$ and that $f$ is
smooth up to $\partial D$ and is extended into the exterior of $D$ constant in
the normal direction.  We next perform an asymptotic analysis to estimate the
rate of convergence of the corresponding approximations.

\subsection{Asymptotic analysis}
\label{sec:asymptotic_analysis}

\begin{figure}[b!p]
  \centering
  \begin{tikzpicture}
    [
    yscale=0.8,
    interface/.style={thick},
    inner/.style={fill=gray,dotted,fill opacity=0.2},
    outer/.style={fill=gray,dashed,fill opacity=0.3},
    labels/.style={above right, font=\small},
    ]

    \begin{scope}[scale=0.8]
      \draw[thick] (0,0) rectangle(9,5);
      \draw[thick] (2,1)
             .. controls (1,1) and (1,4) .. (3,3.5)
             .. controls (5,3) and (4,5) .. (6,4)
             .. controls (8,3) and (8,1.5) .. (7,1.5) node (g1) {}
             .. controls (3,1.5) and (3,1) .. (2,1);
      \node at (7.4,0.4) {$\Omega$};
      \node at (2.0,2.3) {$D$};
      \coordinate (a) at (2.4,1.6);
      \coordinate (b) at (3.4,1.6);
      \coordinate (c) at (2.4,0.6);
      \coordinate (d) at (3.4,0.6);
      \draw[very thin] (c) rectangle (b);
    \end{scope}

    \begin{scope}[xshift=2.5cm,yshift=4.4cm]
      \draw[very thin] (a) -- (0, 3);
      \draw[very thin] (b) -- (8, 3);
      \draw[very thin] (c) -- (0,-3);
      \draw[very thin] (d) -- (8,-3);
      \fill[white] (0,-3) rectangle (8,3);
      \draw[very thin, densely dotted] (a) -- (0, 3);
      \draw[very thin, densely dotted] (b) -- (8, 3);

      \draw[interface] (0,0)
        .. controls (2, 0.5) and (3, 0.5) .. (4,0)
        .. controls (5,-0.5) and (6,-1.0) .. (8,0);

      \draw[inner] (0,2.0)
        .. controls (2,2.5) and (3,2.5) .. (4,2.0)
        .. controls (5,1.5) and (6,1.0) .. (8,2.0) -- (8,-2.0)
        .. controls (6,-3.0) and (5,-2.5) .. (4,-2.0)
        .. controls (3,-1.5) and (2,-1.5) .. (0,-2.0) -- cycle;

      \draw[outer] (0,3.0) -- (8,3.0) -- (8,1.0)
        .. controls (6,0.0) and (5,0.5) .. (4,1.0)
        .. controls (3,1.5) and (2,1.5) .. (0,1.0) -- cycle;
      \draw[outer] (0,-3.0) -- (8,-3.0) -- (8,-1.0)
        .. controls (6,-2.0) and (5,-1.5) .. (4,-1.0)
        .. controls (3,-0.5) and (2,-0.5) .. (0,-1.0) -- cycle;

      \node[labels] at (0.1,2.25) {Outer region};
      \node[labels] at (0.1,1.25) {Overlapping region};
      \node[labels] at (0.1,0.25) {Inner region};
      \draw[decorate,decoration=brace] (8.1, 3.0) --
        node[right=0.5em] {$D$} (8.1, 0.1);
      \draw[decorate,decoration=brace] (8.1,-0.1) --
        node[right=0.5em] {$\Omega$} (8.1,-3.0);
      \node[right=0.5em] at (8.1,0) {$\partial D$};
    \end{scope}
  \end{tikzpicture}
  \caption{A sketch of the regions used for the matched asymptotic expansions.
    The inner region is marked with a light grey colour and the outer region
    with a slightly darker grey colour.  The overlapping region is marked with
    the darkest grey colour.}
  \label{fig:regions}
\end{figure}
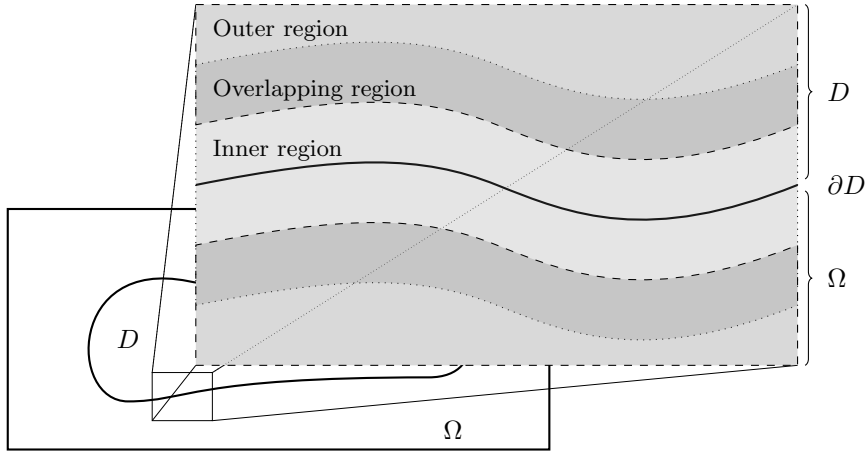

To show asymptotic convergence, we need to consider the expansions of the
diffuse-domain variables in powers of the interface thickness $\epsilon$ in
regions close to and far from the interface.  These are called inner and outer
expansions, respectively.  The two expansions are then matched in a region
where both are valid, see \cref{fig:regions}, which provides the boundary
conditions for the outer variables.  We refer the reader to \cite{Caginalp88}
and \cite{Pego88} for more details and discussion of the general procedure.

The outer expansion for the variable $u(\vct x;\epsilon)$ is simply
\begin{equation}
  u(\vct x;\epsilon) = u^{(0)}(\vct x)
    + \epsilon   u^{(1)}(\vct x)
    + \epsilon^2 u^{(2)}(\vct x) + \cdots.
  \label{eq:outer}
\end{equation}
The outer expansion of an equation is then found by inserting the expanded
variables into the equation.

The inner expansion is found by introducing a local coordinate system near the
interface $\partial D$,
\[
  \vct x(\vct s,z;\epsilon) = \vct X(\vct s;\epsilon)
    + \epsilon z\vct n(\vct s;\epsilon),
\]
where $\vct X(\vct s;\epsilon)$ is a parametrization of the interface, $\vct
n(\vct s;\epsilon)$ is the interface normal vector that points out of $D$, $z$
is the stretched variable
\[
  z = \frac{r(\vct x)}{\epsilon},
\]
and $r$ is the signed distance from the point $\vct x$ to $\partial D$.  In the
local coordinate system, the derivatives become
\begin{align*}
  \grad &= \frac{1}{\epsilon}\vct n\partial_z
           + \frac{1}{1+\epsilon z\kappa}\grads, \\
  \lapl &= \frac{1}{\epsilon^2}\partial_{zz}
           + \frac{1}{\epsilon}\frac{\kappa}{1+\epsilon z\kappa}\partial_z
           + \frac{1}{1+\epsilon z\kappa}\grads
             \cdot \left(\frac{1}{1+\epsilon z\kappa}\grads\right),
\end{align*}
where $\kappa\equiv\divs\vct n$ is the curvature of the interface.  Note that
$\vct n = -\frac{\nabla\phi}{|\nabla\phi|}$.  The inner variable $\hat u(z,\vct
s;\epsilon)$ is now given by
\[
  \hat u(z,\vct s;\epsilon) \equiv u(\vct x;\epsilon)
  = u(\vct X(\vct s;\epsilon) + \epsilon z\vct n(\vct s;\epsilon);\epsilon),
\]
and the inner expansion is
\begin{equation}
  \hat u(z,\vct s;\epsilon) = \hat u^{(0)}(z,\vct s)
    + \epsilon \hat u^{(1)}(z,\vct s)
    + \epsilon^2\hat u^{(2)}(z,\vct s) + \cdots.
  \label{eq:inner}
\end{equation}

To obtain the matching conditions, we assume that there is a region of overlap
where both the expansions are valid.  In this region, the solutions have to
match.  In particular, if we evaluate the outer expansion in the inner
coordinates, this must match the limits of the inner solutions away from the
interface, that is
\[
  u(\vct X + \epsilon z\vct n;\epsilon) \simeq \hat u(z,\vct s;\epsilon).
\]
Insert the expansions into \cref{eq:outer,eq:inner} to get
\[
  \begin{split}
  u^{(0)}(\vct X + \epsilon z\vct n)
    + \epsilon u^{(1)}(\vct X + \epsilon z\vct n)
    + \epsilon^2u^{(2)}(\vct X + \epsilon z\vct n) + \cdots \\
  \simeq \hat u^{(0)}(z,\vct s)
    + \epsilon \hat u^{(1)}(z,\vct s)
    + \epsilon^2\hat u^{(2)}(z,\vct s) + \cdots.
  \end{split}
\]
The terms on the left-hand side can be expanded as a Taylor series,
\[
  u^{(k)}(\vct X + \epsilon z\vct n) = u^{(k)}(\vct s)
    + \epsilon z\ndot\grad u^{(k)}(\vct s)
    + \frac{\epsilon^2z^2}{2} \ndot\grad\grad u^{(k)}(\vct s)\cdot\vct n
    + \cdots,
\]
where $k \in \mathbb N$ and $u^{(k)}(\vct s)\equiv u^{(k)}(\vct X(\vct
s;\epsilon))$.  Now we end up with the matching equation
\[
  \begin{split}
    u^{(0)}(\vct s)
      +& \epsilon\left(u^{(1)}(\vct s) + z\ndot\grad u^{(0)}(\vct s)\right) \\
      +& \epsilon^2\left(
        u^{(2)}(\vct s) + z\ndot\grad u^{(1)}(\vct s)
        + \frac{z^2}{2} \ndot\grad\grad u^{(0)}(\vct s)\cdot\vct n \right) \\
      +& \cdots \simeq \hat u^{(0)}(z,\vct s) + \epsilon \hat u^{(1)}(z,\vct s)
      + \epsilon^2\hat u^{(2)}(z,\vct s) + \cdots,
  \end{split}
\]
which must hold when the interface width is decreased, that is $\epsilon \to
0$.  In the matching region it is required that $\epsilon z = \bigo 1$.  Under
this condition, if we let $z \to \pm\infty$, we get the following asymptotic
matching conditions:
\[
  \zlimpm \hat u^{(0)}(z,\vct s) = u^{(0)}(\vct s),
  \label{eq:match1}
\]
and as $z \to \pm\infty$,
\begin{align}
  \label{eq:match2}
  \hat u^{(1)}(z,\vct s) &= u^{(1)}(\vct s)
    + z\ndot\grad u^{(0)}(\vct s) + \smallo 1, \\
  \nonumber
  \hat u^{(2)}(z,\vct s) &= u^{(2)}(\vct s)
    + z\ndot\grad u^{(1)}(\vct s) \\ &\quad
  \label{eq:match3}
    + \frac{z^2}{2} \ndot\grad\grad u^{(0)}(\vct s)\cdot\vct n
    + \smallo 1,
\end{align}
where the quantities on the right-hand side are the limits from the interior
($-$) and exterior ($+$) of $D$.  Here $\smallo 1$ means that the expressions
approach equality when $z\to\pm\infty$.  That is, $\smallo 1$ is defined such
that if some function $f(z) = \smallo 1$, then we have $\zlimpm f(z) = 0$.

\subsection{Poisson equation with Robin boundary conditions}

Now we are ready to consider the Poisson equation with Robin boundary
conditions,
\begin{equation}
  \begin{alignedat}{2}
    \lapl u      &= f               & \txin D, \\
    \ndot\grad u &= k(u - g) \qquad & \txon \partial D,
  \end{alignedat}
  \label{eq:poiss_robin}
\end{equation}
where $k\le 0$.  Consider a general DDM approximation,
\begin{equation}
  \div\left(\phi\grad u\right) + \frac{1}{\epsilon^2}\psi = \phi f.
  \label{eq:general_dda}
\end{equation}
where $\psi$ represents the BC approximation in the DDM.  The scaling factor
$1/\epsilon^2$ is taken for later convenience.  If we assume that $\psi$ is
local to the interface (e.g., vanishes to all orders in $\epsilon$ away from
$\partial D$) and that $f$ is independent of $\epsilon$ (e.g., is smooth in
a neighbourhood of $\partial D$ and is extended constant in the normal
direction out of $D$), which is the case for the approximations BC1 and BC2
given in \cref{eq:bc_robin}, then the outer solution to this equation
when $z\to -\infty$ satisfies
\begin{equation}
  \begin{split}
    \lapl u^{(0)} &= f, \\
    \lapl u^{(1)} &= 0, \\
    \lapl u^{(k)} &= 0,\qquad k = 2,3,\dots.
  \end{split}
  \label{eq:dda_outer}
\end{equation}
Now, if $u^{(0)}$ satisfies \cref{eq:poiss_robin} and $u^{(1)} \ne 0$ then the
outer expansion $u\approx u^{(0)} + \epsilon u^{(1)} +~\dots$ and the DDM is
asymptotically first-order accurate.  However, if $u^{(1)} = 0$, then $u\approx
u^{(0)} + \epsilon^2 u^{(2)} + \dots$ and the DDM is asymptotically
second-order accurate.  Determining which of these is the case requires
matching the outer solutions to the solutions of the inner equations.

\subsubsection{Matching conditions}

Before we analyse the inner expansions, we develop a higher-order matching
condition based on \cref{eq:match2,eq:match3} that matches a Robin boundary
condition for $u^{(1)}$.  First we take the derivative of \cref{eq:match3} with
respect to $z$ and subtract $k$ times \cref{eq:match2}, which gives
\[
  \hat u^{(2)}_z - k\hat u^{(1)}
  = - ku^{(1)} - kz\ndot\grad u^{(0)}
    + \ndot\grad u^{(1)} + z \ndot\grad\grad u^{(0)}\cdot\vct n.
\]
Move the terms that make up a Robin condition for $u^{(1)}$ to the left-hand
side, and move the rest to the right-hand side, that is
\begin{equation}
  \ndot\grad u^{(1)} - ku^{(1)}
  = \hat u^{(2)}_z
    - k\hat u^{(1)}
    + kz\ndot\grad u^{(0)}
    - z \ndot\grad\grad u^{(0)}\cdot\vct n.
  \label{eq:robin_mc}
\end{equation}
The Laplacian can be decomposed into normal and tangential components as
\begin{equation}
  \lapl u = \ndot\grad\grad u\cdot\vct n + \kappa\ndot\grad u + \lapls u,
  \label{laplace decomp}
\end{equation}
which can be shown by writing the gradient vector as $\grad = \vct n\ndot\grad
+ \grads$.  We can therefore write
\[
  \ndot\grad\grad u^{(0)}\dotn
  = f - \kappa\ndot\grad u^{(0)} - \lapls u^{(0)}
  = \hat f^{(0)} - \kappa k\left(u^{(0)} - \hat g\right) - \lapls u^{(0)},
\]
where we have assumed that $u^{(0)}$ satisfies the system
\eqref{eq:poiss_robin}, as demonstrated below.  If we insert this into the
matching condition \eqref{eq:robin_mc}, we get
\begin{equation}
  \ndot\grad u^{(1)} - ku^{(1)}
  = \hat u^{(2)}_z
    - k\hat u^{(1)}
    - z\left(\hat f^{(0)}
      - \left(\kappa+k\right)k\left(u^{0)}-\hat g\right)
      - \lapls u^{(0)}
      \right),
  \label{eq:robin_mc2}
\end{equation}
as $z\to - \infty$.

\subsubsection{Inner expansions}
\label{sec:dda_robin_inner}

Now consider the inner expansion of \cref{eq:general_dda},
\[
  \frac{1}{\epsilon^2}\left(\phi\hat u_z\right)_z
    + \frac 1 \epsilon\frac{\kappa}{1+\epsilon z\kappa}\phi\hat u_z
    + \frac{\phi}{1+\epsilon z\kappa}\grads
      \cdot\left(\frac{1}{1+\epsilon z\kappa}\grads\hat u\right)
    + \frac{1}{\epsilon^2}\hat\psi
  = \phi\hat f.
\]
Expand $\hat u$, $\hat f$, $\hat \psi$ and $\displaystyle{\frac{1}{1+\epsilon
    z\kappa}}$ in powers of $\epsilon$, to get
\begin{multline*}
  \frac{1}{\epsilon^2}\left(\phi\hat u_z^{(0)}\right)_z
    + \frac{1}{\epsilon}\left(\phi\hat u_z^{(1)}\right)_z
    + \left(\phi\hat u_z^{(2)}\right)_z
    + \frac 1 \epsilon\kappa\phi\hat u_z^{(0)}
    + \kappa\phi\hat u_z^{(1)}-z\kappa^2\hat u_z^{(0)} \\
    + \phi\lapls\hat u^{(0)}
    + \frac{1}{\epsilon^2}\hat\psi^{(0)}
    + \frac{1}{\epsilon}\hat\psi^{(1)}
    + \hat\psi^{(2)}
  = \phi\hat f^{(0)} + \bigo{\epsilon}.
\end{multline*}
and then collect the leading order terms.  Note that because $f$ is smooth up
to $\partial D$ and extended constantly outside $D$ we have that $\hat f^{(0)}$
is independent of $z$.  The lowest power of $\epsilon$ gives
\[
  \left( \phi \hat u^{(0)}_z \right)_z = - \hat\psi^{(0)}.
\]
Suppose that $\hat\psi^{(0)} = 0$, which is the case as we show below for BC1
and BC2, then we obtain $\hat u_z^{(0)} = 0$.  By the matching condition
\eqref{eq:match1}, this gives $\hat u^{(0)}(z,\vct s) = u^{(0)}(\vct s)$,
where $u^{(0)}(\vct s)$ is the limiting value of $u^{(0)}$.

The next order terms give
\begin{equation}
  \left( \phi\hat u_z^{(1)} \right)_z = - \hat\psi^{(1)}.
  \label{eq:first_bc}
\end{equation}
Integrating from $-\infty$ to $+\infty$ in $z$ and using the matching condition
\eqref{eq:match2}, we get
\[
  \ndot\grad u^{(0)} = \eint{\hat\psi^{(1)}}.
\]
To obtain a Robin boundary condition for $u^{(0)}$, we need that
\[
  \eint{\hat\psi^{(1)}} = k(u^{(0)} - g).
\]

Now consider the zeroth order terms,
\begin{equation}
  \left( \phi\hat u_z^{(2)} \right)_z
  = \phi\hat f^{(0)}
    - \kappa\phi\hat u_z^{(1)}
    - \phi\lapls u^{(0)}
    - \hat\psi^{(2)}.
  \label{eq:ddarobin0}
\end{equation}
If we subtract
\[
  \left(\phi k\hat u^{(1)} + z\phi\left(
    \hat f^{(0)} - (\kappa+ k)k\left(u^{(0)}-\hat g\right) - \lapls u^{(0)}
  \right)\right)_z
\]
from both sides of \cref{eq:ddarobin0}, we get
\begin{multline*}
  \left(\phi\hat u^{(2)}_z
    - \phi k\hat u^{(1)}
    - z\phi\left(\hat f^{(0)}
      - (\kappa+k) k\left(u^{(0)}-\hat g\right)
      - \lapls u^{(0)}
     \right) \right)_z \\
  = - \hat\psi^{(2)} - k\phi_z\hat u^{(1)}
    - z\phi_z\left(\hat f^{(0)}
      - (\kappa +k)k\left(u^{(0)}-\hat g\right)
      - \lapls u^{(0)} \right) \\
    - \phi(\kappa + k)\left(\hat u^{(1)}_z-k\left(u^{(0)}-g\right)\right),
\end{multline*}
where we have taken into account the cancellation of terms and used the fact
that $\hat f^{(0)}$ and $\hat g$ are independent of $z$.  The latter
holds when $g$ is extended as a constant in the normal direction off
$\partial D$, e.g., $\hat g(z,s) = g(s)$ and is independent of $z$ and
$\epsilon$.  Next, we integrate and use the matching condition
\eqref{eq:robin_mc2} on the left-hand side,
\begin{multline}
  \ndot\grad u^{(1)} - k u^{(1)}
  = \eint{\left(\hat\psi^{(2)} + k\phi_z\hat u^{(1)}
      + \phi(\kappa + k)\left(\hat u^{(1)}_z-k\left(u^{(0)}-g\right)\right)
      \right)} \\
    + \left(\hat f^{(0)} - (\kappa + k)k\left(u^{(0)} - \hat g\right)
      - \lapls u^{(0)}\right) \eint{z\phi_z}.
  \label{order 1 BC}
\end{multline}
If the right-hand side of \cref{order 1 BC} vanishes, then we obtain
$\ndot\grad u^{(1)} - k u^{(1)}=0$ from which we can conclude that the outer
solution $u^{(1)} = 0$ since $u^{(1)}$ is harmonic: $\lapl u^{(1)} = 0$.  Next
we analyse the boundary condition approximations BC1 and BC2.

\subsubsection{Analysis of BC1}

The BC1 approximation corresponds to
\[
  \frac{1}{\epsilon^2}\psi = k(u-g)|\nabla\phi|.
\]
Since $\phi=1$ in the outer region (interior part of $D$), we conclude that
$\psi$ vanishes in the outer region.  In the inner region, we have
\begin{equation}
  \frac{1}{\epsilon^2}\psi
  = -\frac{k}{\epsilon}\left(\hat u - g\right)\phi_z
  = -\frac{k}{\epsilon}\left(\hat u^{(0)}- g\right)\phi_z
     - k\hat u^{(1)}\phi_z + \bigo{\epsilon},
  \label{expansion of BC1}
\end{equation}
where we have used that $\hat g(z,s) = g(s)$ and is independent of $z$ and
$\epsilon$.  Since $\hat\psi^{(0)} = 0$, we conclude from our analysis in
\cref{sec:dda_robin_inner} that $\hat u^{(0)}_z = 0$ and hence $\hat u^{(0)}
= u^{(0)}$.

The next orders of \cref{expansion of BC1} give
\[
  \hat\psi^{(1)} = -k\left(u^{(0)} - g\right)\phi_z,
  \qquad\text{and}\quad
  \hat\psi^{(2)} = -k\hat u^{(1)}\phi_z.
\]
A direct calculation then shows that
\[
  \eint{\hat\psi^{(1)}} = k(u^{(0)} - g),
\]
as desired.  Thus, the leading order outer solution $u^{(0)}$ satisfies the
problem (\ref{eq:poiss_robin}).

To continue, we must first consider \cref{eq:first_bc},
\[
  \left(\phi\hat u_z^{(1)}\right)_z = k\left(u^{(0)} - g\right)\phi_z,
\]
from which we get
\[
  \hat u^{(1)}(z,\vct s) = u^{(1)}(\vct s) + zk\left(u^{(0)} - g\right),
\]
where $u^{(1)}(\vct s)$ is the limiting value of the outer solution (e.g., see
\cref{eq:match2}).  Combining this with \cref{order 1 BC}, we obtain
\begin{equation}
  \ndot\grad u^{(1)} - k u^{(1)} = \left(\hat f^{(0)}
      - (\kappa + k)k\left(u^{(0)} - g\right)
      - \lapls\hat u^{(0)}\right)\eint{z\phi_z}.
  \label{order 1 BC 1}
\end{equation}
Further, it follows from the definition of the phase-field
function~\eqref{eq:characteristic} that $z\phi_z$ is an odd function.
Therefore the integral on the right-hand side of \cref{order 1 BC 1} is
equal to zero.  Thus $\ndot\grad u^{(1)}-k u^{(1)} = 0$ and so by our
arguments below \cref{eq:dda_outer}, the DDM with BC1 is second-order accurate
in $\epsilon$.

\subsubsection{Analysis of BC2}

When the BC2 approximation is used, we obtain
\[
  \frac{1}{\epsilon^2}\psi = \epsilon k(u-g)|\nabla\phi|^2.
\]
Accordingly, in the inner region, we obtain
\[
  \hat\psi^{(0)} = 0,
  \qquad
  \hat\psi^{(1)} = k\left( u^{(0)} - g\right)\phi_z^2,
  \qquad\text{and}\quad
  \hat\psi^{(2)} = k \hat u^{(1)}\phi_z^2.
\]
Since $\eint{\phi_z^2} = 1$, we get
\[
 \eint{\hat\psi^{(1)}} = k(u^{(0)} - \hat g),
\]
as desired.  From \cref{eq:first_bc} we have
\[
  \left(\phi\hat u_z^{(1)}\right)_z = -k\left(u^{(0)} - g\right)\phi_z^2.
\]
Using that $\phi_z = -6\phi(1-\phi)$, this gives
\begin{equation}
  \hat u^{(1)}(z,\vct s) = C(\vct s) - k\left(u^{(0)} - g\right) F(\phi),
  \label{u1 inner solve B2}
\end{equation}
where
\begin{align}
  \label{u1 inner solve B2 a}
  F(\phi) &= -\frac{1}{6}\log\left(1 - \phi\right) + \frac{\phi}{3}, \\
  \nonumber
  C(\vct s) &= u^{(1)}(\vct s) + \frac{k}{3}\left(u^{(0)} - g\right).
\end{align}
Combining \cref{u1 inner solve B2,u1 inner solve B2 a,order 1 BC} we get
\begin{multline}
  \ndot\grad u^{(1)} - k u^{(1)}
    = kC(\vct s)\eint{\left(\phi_z^2 + \phi_z\right)} \\
    - k^2\left(u^{(0)} - \hat g\right)
      \eint{F(\phi)\left(\phi_z^2 + \phi_z\right)} \\
    + \left(\kappa + k\right)k\left(u^{(0)} - \hat g\right)
      \eint{\phi\left(3\phi - 2\phi^2 - 1\right)}.
  \label{order 1 BC 2}
\end{multline}
Direct calculations show that
\[
  \eint{\left(\phi_z^2 + \phi_z\right)}
  = \eint{\left(3\phi^2 - 2\phi^3 - \phi\right)} = 0,
\]
and
\[
  \eint{F(\phi)\left(\phi_z^2 + \phi_z\right)} = -\frac{1}{36}.
\]
Using these in \cref{order 1 BC 2}, we get
\begin{equation}
  \ndot\grad u^{(1)} - k u^{(1)} = \frac{1}{36}k^2\left(u^{(0)} - g\right).
  \label{order 1 BC 2 final}
\end{equation}
This shows that the DDM with BC2 is only first-order accurate because the
solution $u^{(1)}$ of
\[
  \begin{alignedat}{2}
    \lapl u^{(1)}                 &= 0 & \txin D, \\
    \ndot\grad u^{(1)} - ku^{(1)}
      &= \frac{1}{36}k^2\left(u^{(0)} - g\right) \qquad
      & \txon \partial D,
  \end{alignedat}
\]
is in general not equal to $0$, e.g., $u^{(1)}\ne 0$.

\subsubsection{Analysis of a second-order modification of BC2}
\label{sec:dda_robin_bc2m}

In order to modify BC2 to achieve second-order accuracy, we introduce
$\tilde\psi$ such that
\[
  \hat \psi^{(0)} = 0,
  \qquad
  \hat \psi^{(1)} = k\left( u^{(0)}- g\right)\phi_z^2,
  \qquad\text{and}\quad
  \hat \psi^{(2)} = k \hat u^{(1)}\phi_z^2 +\hat{\tilde\psi}^{(0)}.
\]
That is, $\tilde\psi$ perturbs only the higher order terms in the inner
expansion and is chosen to cancel the term on the right-hand side of
\cref{order 1 BC 2 final} in order to achieve $\ndot\grad u^{(1)} - k u^{(1)}
= 0$, which in turn implies that $u^{(1)} = 0$ and the new formulation is
second-order accurate.  The correction $\tilde\psi$ does not affect the
$O(\epsilon^{-2})$ or $O(\epsilon^{-1})$ orders in the system.  Thus, $\hat
u^{(0)}$ and $\hat u^{(1)}$ are unchanged from the previous subsection.
\Cref{order 1 BC} now becomes
\[
 \ndot\grad u^{(1)} - k u^{(1)}
   = \frac{1}{36}k^2\left(u^{(0)} - g\right)
     + \eint{\hat{\tilde\psi}^{(0)}},
\]
so we wish to determine $\hat{\tilde\psi}^{(0)}$ such that
\[
  \eint{\hat{\tilde\psi}^{(0)}}
  = -\frac{1}{36}k^2\left(u^{(0)} - g\right).
\]
Two simple ways of achieving this are to take
\[
  \hat{\tilde\psi}^{(0)} = -\frac{1}{36}k^2\left(u^{(0)} - g\right)
    \times
    \begin{cases}
      -\phi_z,
      \quad\text{or} \\
      \phi_z^2.
    \end{cases}
\]
Putting everything together, we can obtain BC2M, a second-order version of BC2,
using
\[
  \text{BC2M} =
    \begin{cases}
      \text{BC2M1} = \epsilon k(u - g) |\nabla\phi|
        \left(|\nabla\phi| - \frac{k}{36}\right),
      \quad\text{or} \\
      \text{BC2M2} = \epsilon k(u - g) |\nabla\phi|^2
        \left(1 - \epsilon\frac{k}{36}\right).
    \end{cases}
\]
The resulting DDM is an elliptic system since $k<0$, as required
for the Robin boundary condition.  In each
instance, this is guaranteed if the interface thickness $\epsilon$ is
sufficiently small.

\subsubsection{Other approaches to second-order BCs}
\label{sec:dda_robin_other}

Thus far, we have taken advantage of integration to achieve second-order
accuracy.  Alternatively, one may try to add correction terms to directly
obtain second-order boundary conditions without relying on integration.  For
example, from \cref{order 1 BC} to achieve second-order accuracy we may take
\begin{multline}
  \hat\psi^{(2)}
  = -k\phi_z\hat u^{(1)} - \phi(\kappa+k)
      \left(\hat u^{(1)}_z - k\left(u^{(0)} - g\right)\right) \\
    - \left(\hat f^{(0)}
      - (\kappa + k)k\left(u^{(0)} - \hat g\right)
      - \lapls\hat u^{(0)}\right) z\phi_z,
  \label{local correction}
\end{multline}
where $\hat u^{(1)}$ is a functional of $\hat\psi^{(1)}$.  This provides
another prescription of how to obtain a second-order accurate boundary
condition, which could in principle lead to faster asymptotic convergence since
it directly cancels a term in the inner expansion of the asymptotic matching.
As an illustration, let us use BC1 as a starting point even though this
boundary condition is already second-order accurate.  Through the prescription
in \cref{local correction} above, we derive another second-order accurate
boundary condition.  To see this, write
\[
  \hat\psi^{(0)} = 0,
  \qquad
  \hat \psi^{(1)} = -k\left(u^{(0)} - g\right)\phi_z,
  \qquad\text{and}\quad
  \hat \psi^{(2)} = -k\hat u^{(1)}\phi_z + \hat{\tilde\psi}^{(0)},
\]
then from \cref{local correction} we get
\[
  \hat{\tilde\psi}^{(0)} = -\left(\hat f^{(0)}
    - (\kappa +k)k\left(u^{(0)} - \hat g\right)
  - \lapls\hat u^{(0)}\right) z\phi_z.
\]
This can be achieved by taking
\[
  \frac{1}{\epsilon^2}\psi = k(u-g) |\nabla\phi| \left(1 - r(k + \kappa)\right)
    + r|\nabla\phi|\left(f - \lapls u\right).
\]
where $r$ is the signed distance to $\partial D$ as defined earlier.  Note that
we can also achieve second-order accuracy by taking instead
\[
  \hat\psi^{(2)} = -k\phi_z\hat u^{(1)}
  - \phi(\kappa + k)\left(\hat u^{(1)}_z - k\left(u^{(0)} - g\right)\right)
  - \left(\hat f^{(0)} - (\kappa + k)k\left(u^{(0)} - \hat g\right)\right)
      z\phi_z,
\]
where we use the fact that the integral involving $\lapls u^{(0)}$ vanishes in
\cref{order 1 BC}.  We refer to these choices, which are by no means
exhaustive, as
\begin{equation}
  \text{BC1M} =
    \begin{cases}
      \text{BC1M1} = k(u-g) |\nabla\phi| \left(1 - r(k + \kappa)\right)
        + r|\nabla\phi| \left(f - \lapls u\right),
      \quad\text{or} \\
      \text{BC1M2} = k(u-g) |\nabla\phi| \left(1 - r(k + \kappa)\right)
        + r|\nabla\phi| f.
    \end{cases}
  \label{BC1M versions}
\end{equation}
We remark, however, that this prescription may not always lead to an optimal
numerical method.  For example, when using BC1M2, the system is guaranteed to
be elliptic when $1 - r(k + \kappa) > 0$ for $|r|\approx \epsilon$, which
puts an effective restriction on the interface thickness $\epsilon$ depending
on the values of $k$ and $\kappa$.  When BC1M1 is used, the situation is more
delicate since ellipticity cannot be guaranteed when $r>0$ due to the $\lapls
u$ term.  Recall that $r>0$ outside the original domain $D$ and so this issue
is associated with the extending the modified boundary condition outside $D$.
In future work, we plan to consider different extensions that automatically
guarantee ellipticity.

To summarize, we have shown that the DDM
\[
  \div\left(\phi\grad u\right) + \text{BC} = \phi f
\]
is a second-order accurate approximation of the system \eqref{eq:poiss_robin}
when BC1, BC2M, and BC1M are used.  When BC2 is used, the DDM is only
first-order accurate.

\subsection{Reaction-diffusion equation with Neumann boundary conditions}

Since the Poisson equation with Neumann boundary conditions does not have
a unique solution, we instead consider the steady reaction-diffusion equation
with Neumann boundary conditions,
\begin{equation}
  \begin{alignedat}{2}
    \lapl u - u  &= f        & \txin D, \\
    \ndot\grad u &= g \qquad & \txon \partial D.
  \end{alignedat}
  \label{eq:poiss_neumann}
\end{equation}

Again we consider a general DDM approximation,
\begin{equation}
  \div\left(\phi\grad u\right) - \phi u + \frac{1}{\epsilon^2}\psi = \phi f.
  \label{eq:reacdiff_general_dda}
\end{equation}
Under the same conditions on $\psi$ as in the previous section, the outer
solution now satisfies
\begin{align*}
  \lapl u^{(0)} - u^{(0)} &= f, \\
  \lapl u^{(k)} - u^{(k)} &= 0,\qquad k = 1,2,3,\dots
\end{align*}
As in the Robin case, if $u^{(0)}$ satisfies \cref{eq:poiss_neumann} and
$u^{(1)} \ne 0$ then the DDM is first-order accurate.  However, if $u^{(1)}
= 0$, then the DDM is second-order accurate.

To construct the boundary condition for $u^{(1)}$, we follow the approach from
the Robin case and combine \cref{eq:match3,laplace decomp} to get
\[
  \ndot\grad\grad u^{(0)}\dotn
    = \hat f^{(0)} - \kappa \hat g - \lapls u^{(0)} + u^{(0)},
\]
assuming that $u^{(0)}$ satisfies the system \eqref{eq:poiss_neumann} as
demonstrated below, and to get
\[
  \ndot\grad u^{(1)} = \hat u_z^{(2)}
    - z\left(\hat f^{(0)} - \kappa\hat g - \lapls u^{(0)} + u^{(0)}\right),
\]
as $z\to -\infty$.

\subsubsection{Inner expansions}

The inner expansion of \cref{eq:reacdiff_general_dda} is analogous to the Robin
case derived in \cref{sec:dda_robin_inner}.  As before, if $\hat\psi^{(0)} = 0$
then $\hat u_z^{(0)}
= 0$ and $\hat u^{(0)}(z,\vct s) = u^{(0)}(\vct s)$, the limiting value of the
outer solution.  \Cref{eq:first_bc} still holds at the next order and so to get
the desired boundary condition for $u^{(0)}$, we need
\begin{equation}
  \eint{\hat\psi^{(1)}} = g.
  \label{integral constraint Neumann}
\end{equation}
Analogously to \cref{eq:ddarobin0} the next order equation is
\begin{equation}
  \left(\phi\hat u^{(2)}_z\right)_z
  + \phi\kappa\hat u_z^{(1)}
  + \phi\lapls \hat u^{(0)}
  - \phi\hat u^{(0)}
  + \hat\psi^{(2)}
  = \phi\hat f^{(0)}.
  \label{eq:o1}
\end{equation}
Subtracting
\[
  -\left(z\phi\left(\hat f^{(0)}
         - \kappa\hat g - \lapls u^{(0)} + u^{(0)}\right)\right)_z
\]
from \cref{eq:o1} we get
\[
  \begin{split}
    \bigg(\phi&\hat u_z^{(2)}
      - z\phi\left(\hat f^{(0)} - \kappa\hat g
      - \lapls u^{(0)} + u^{(0)}\right)\bigg)_z \\
    &= - \hat\psi^{(2)} - \phi\kappa\left(\hat u^{(1)}_z - \hat g\right)
      - z\phi_z\left(\hat f^{(0)} - \kappa \hat g
        - \lapls u^{(0)} + u^{(0)}\right),
  \end{split}
\]
where we have used $\hat u^{(0)}=u^{(0)}$ as justified below.  Integrating, we
obtain
\begin{equation}
  \ndot\grad u^{(1)} = \eint{\left(
      \hat\psi^{(2)} + \phi\kappa\left(\hat u^{(1)}_z - \hat g\right)
      + z\phi_z\left(\hat f^{(0)}
        - \kappa\hat g
        - \lapls u^{(0)} + u^{(0)}\right)\right)}.
  \label{order 1 BC Neumann}
\end{equation}
As in the Robin case, if the right-hand side of \cref{order 1 BC Neumann}
vanishes then we may conclude that $u^{(1)} = 0$ since $u^{(1)}$ satisfies
$\lapl u^{(1)} - u^{(1)} = 0$ with zero Neumann boundary conditions.  We next
analyse the boundary conditions BC1 and BC2.

\subsubsection{Analysis of BC1}

When the BC1 approximation is used, we obtain
\[
  \frac{1}{\epsilon^2}\psi = g|\nabla\phi|,
\]
and
\begin{equation}
  \hat \psi^{(0)} = 0,
  \qquad
  \hat \psi^{(1)} = -g\phi_z,
  \qquad\text{and}\quad
  \hat \psi^{(2)} = 0.
  \label{psi expansion for BC1 Neumann}
\end{equation}
Accordingly, we find that $\hat u^{(0)}(z,\vct s) = u^{(0)}(\vct s)$ and
\cref{integral constraint Neumann} holds.  Thus $u^{(0)}$ satisfies the system
\eqref{eq:poiss_neumann} as claimed above.

At the next order, from \cref{eq:first_bc,psi expansion for BC1 Neumann} we
obtain
\begin{equation}
  \hat u^{(1)}(z,\vct s) = u^{(1)}(\vct s) + z\hat g.
  \label{hat u1 Neumann}
\end{equation}
Thus, combining \cref{hat u1 Neumann,order 1 BC Neumann} we get
\[
  \ndot\grad u^{(1)} = \left(\hat f^{(0)}
    - \kappa\hat g
    - \lapls\hat u^{(0)}+u^{(0)}\right)\eint{z\phi_z} = 0,
\]
from which we conclude that $u^{(1)} = 0$ and the DDM with BC1 is second-order
accurate.

\subsubsection{Analysis of BC2}

When the BC2 approximation is used, we obtain
\[
  \frac{1}{\epsilon^2}\psi = \epsilon g|\nabla\phi|^2,
\]
and
\begin{equation}
  \hat \psi^{(0)} = 0,
  \qquad
  \hat \psi^{(1)} = g\phi_z^2,
  \qquad\text{and}\quad
  \hat \psi^{(2)} = 0.
  \label{psi expansion for BC2 Neumann}
\end{equation}
Analogously to the case when BC1 is used, $\hat u^{(0)}(z,\vct s)
= u^{(0)}(\vct s)$, \cref{integral constraint Neumann} holds, and $u^{(0)}$
satisfies the system \eqref{eq:poiss_neumann}.  At the next order, from
\cref{eq:first_bc,psi expansion for BC2 Neumann} we obtain
\begin{equation}
  \hat u^{(1)}(z,\vct s) = u^{(1)}(\vct s) + z\hat g\left(3\phi
    - 2\phi^2\right).
  \label{hat u1 Neumann BC2}
\end{equation}
Combining \cref{hat u1 Neumann BC2,order 1 BC Neumann} we get
\[
  \ndot\grad u^{(1)}
  = \kappa\hat g\eint{\left(3\phi^2 - 2\phi^3 - \phi\right)}
    + \left(\hat f^{(0)} - \kappa\hat g \lapls\hat u^{(0)} + u^{(0)}\right)
      \eint{z\phi_z} = 0,
\]
from which we conclude that $u^{(1)} = 0$ and the DDM with BC2 is second-order
accurate for the Neumann problem as well, which is different from the Robin
case.

\subsubsection{Other approaches to second-order BCs}

Analogously to the Robin case, to achieve second-order accuracy we may also
take
\[
 \hat\psi^{(2)} + \phi\kappa\left(\hat u^{(1)}_z - \hat g\right)
  + z\phi_z\left(\hat f^{(0)} - \kappa\hat g \lapls u^{(0)}+u^{(0)}\right) = 0.
\]
Following the same reasoning, alternative boundary conditions analogous to
those in \cref{BC1M versions} may be derived
\[
  \text{BC1M} =
    \begin{cases}
      \text{BC1M1} = g|\nabla\phi| \left(1 - r(k + \kappa)\right)
        + r|\nabla\phi| \left(f - \lapls u\right),
      \quad\text{or} \\
      \text{BC1M2} = g|\nabla\phi| \left(1 - r(k + \kappa)\right)
        + r|\nabla\phi| f.
    \end{cases}
\]
Note that as in the Robin case, when BC1M1 is used ellipticity cannot be
guaranteed when $r>0$ due to the $\lapls u$ term.

To summarize, we have shown that the DDM
\[
  \div\left(\phi\grad u\right) - \phi u + \text{BC} = \phi f
\]
is a second-order accurate approximation of the system \eqref{eq:poiss_robin}
when BC1, BC2, and BC1M are used.

\section{Discretizations and numerical methods}
\label{sec:discretization}

The equations are discretized on a uniform grid with the second-order
central-difference scheme.  The discrete system is solved using a multigrid
method, where a red-black Gauss-Seidel type iterative method is used to relax
the solutions (see \cite{Wise07}).  The equations are solved in two-dimensions
in a domain $\Omega = [-2,2]^2$ for all the test cases.  Periodic boundary
conditions are used on the domain boundaries $\partial\Omega$.  The iterations
are considered to be converged when the residual of the current solution has
reached a tolerance of $10^{-9}$.

Since the phase-field function quickly tends to zero outside the physical
domain $D$, it must be regularized in order to prevent the equations from
becoming ill-posed.  We therefore use the modified phase-field function
\[
  \hat\phi = \tau + (1 - \tau)\phi,
\]
where the regularization parameter is set to $\tau = 10^{-6}$ unless otherwise
specified.  In addition, one should note that the chosen boundary condition for
the computational domain, $\Omega$, should not interfere with the physical
domain.  Thus one has to make sure that the distance from the computational
wall to the diffuse interface of $D$ is large enough not to affect the results.

As shown in \cref{sec:asymptotic_analysis}, the normal vector and the
curvature can be calculated from the phase-field function as
\[
  \vct n = -\frac{\grad\phi}{|\grad\phi|},
\]
and
\[
  \kappa = -\div\frac{\grad\phi}{|\grad\phi|}.
\]
The surface Laplacian can be found from the identity
\[
  \lapls \equiv
  \left(I - \vct n\vct n\right) \div\left(I - \vct n\vct n\right)\grad,
\]
where
\[
  (I - \vct n\vct n)\grad \equiv (\delta_{ij} - n_in_j)\partial x_i.
\]
In 2D we get
\begin{align*}
  \nonumber \lapls u
  &= \left(n_1n_2(n_1n_2)_x + n_1n_2(n_1^2)_y - (1-n_1^2)(n_1^2)_x
           - (1-n_2^2)(n_1n_2)_y\right)u_x \\
  &\qquad + \left(n_1n_2(n_1n_2)_y + n_1n_2(n_2^2)_x - (1-n_2^2)(n_2^2)_y
                  - (1-n_1^2)(n_1n_2)_x\right)u_y \\
  &\qquad + \left(\left(1 - n_1^2 \right)^2 + n_1^2n_2^2\right)u_{xx} \\
  &\qquad + \left(\left(1 - n_2^2 \right)^2 + n_1^2n_2^2\right)u_{yy} \\
  &\qquad - 2n_1n_2\left(\left(1 - n_1^2 \right)
                         + \left(1 - n_2^2\right)\right)u_{xy}.
\end{align*}

Below, we verify the accuracy of our numerical implementation on several test
problems in which we manufacture a solution to the DDM with different choices
of boundary conditions through particular choices of $f$.  As suggested by our
analysis, we find that when we include the surface Laplacian, we are unable to
solve the discrete system using the multigrid method even though the correction
term and the subsequent loss of ellipticity outside $D$ is confined to the
interfacial region.  As also mentioned in \cref{sec:dda_robin_other}, future
work involves developing
alternative extensions of the boundary conditions outside $D$ that maintain
ellipticity.  Nevertheless, as a proof of principle, we still consider the
effect of this term by using the surface Laplacian of the analytic solution in
the DDM equations.

\section{Results}
\label{sec:results}

We next investigate the performance of the DDM with different choices of
boundary conditions and compare the results with the exact solution of the
sharp-interface equations for the reaction-diffusion equation with Neumann
boundary conditions and the Poisson equation with Robin boundary conditions.
We consider four different cases with Neumann boundary conditions and three
different cases with Robin boundary conditions.  For each case, we calculate
and compare the error between the calculated solution $u$ and an analytic
solution $u_{\text{an}}$ of the original PDE, which is extended from $D$ into
$\Omega$.  The error is defined as
\[
  E_\epsilon = \frac{\|\phi(u_{\text{an}} - u)\|}{\|\phi u_\text{an}\|},
\]
where $\|\cdot\|$ is a norm and $\phi$ is used to restrict the error to the
physical domain $D$.  The convergence rate in $\epsilon$ as $\epsilon\to 0$ is
calculated as
\[
  k = \log\left(\frac{E_{\epsilon_i}}{E_{\epsilon_{i-1}}}\right) /
  \log\left(\frac{\epsilon_i}{\epsilon_{i-1}}\right),
\]
for a decreasing sequence $\epsilon_i$.  In the following results we mainly use
the $L^2$ norm,
\[
  \|\psi\|_2 = \frac{\sqrt{\sum_{i=1}^N \psi_i^2}}{N},
\]
where $\psi$ is an array with $N$ elements.  For a couple of cases, we also
present the results with the $L^\infty$ norm,
\[
  \|\psi\|_\infty = \max_{i=1}^N\,|\psi_i|.
\]

For a given $\epsilon$, the error $E_\epsilon$ in both $L^2$ and $L^\infty$ is
calculated by refining the grid spacing until a minimum of two leading digits
converge (e.g., stop changing under refinement) in the $L^2$ norm.  In some
cases for the smallest values of $\epsilon$,  the error has not yet converged
to two digits in the $L^\infty$ norm. However, due to memory limits on our
computers, we were not able to use more than $n = 8192$ cells in each
direction.  This limited our ability to obtain grid convergence, particularly
when BC2 (see below) is used for very small values of $\epsilon$.

\subsection{Neumann boundary conditions}

Consider the steady reaction-diffusion equation with Neumann boundary
conditions,
\[
  \begin{alignedat}{3}
    \lapl u - u  &= f       && \txin D, \\
    \ndot\grad u &= g \quad && \txon \partial D.
  \end{alignedat}
\]
In this section we solve the DDM systems
\[
  \div\left(\phi\grad u\right) - \phi u + \text{BC} = \phi f,
\]
where BC refers to selected boundary condition approximations considered in the
previous section.  In the case of BC1M1, as remarked above, the surface
Laplacian term is not solved, rather the surface Laplacian of the analytic
solution is used and is treated as a known source term.

\subsubsection{Case 1}

Consider the case where $D$ is a circle of radius $R=1$ centred at $(0,0)$,
and where the analytic solution to the reaction-diffusion equation in $D$ is
\[
  u_{\text{an}}(x,y) = \frac{1}{4}\left( x^2 + y^2 \right).
\]
This corresponds to $f = 1 - (x^2 + y^2)/4$, $g = 1/2$, and $\lapls
u_{\text{an}} = 0$.  In this case, the curvature is $\kappa = 1$.

\subsubsection{Case 2}

Now consider the case where $D$ is the square $D=[-1,1]^2$.  Again let the
analytic solution in $D$ be
\[
  u_{\text{an}}(x,y) = \frac{1}{4}\left( x^2 + y^2 \right),
\]
so that $f = 1 - (x^2 + y^2)/4$, $g = 1/2$, and $\lapls u_{\text{an}} = 1/2$.
In this case the curvature is zero almost everywhere.

To initialize the square domain $D$, the signed-distance function is defined as
\[
  r(x,y) =
  \begin{cases}
    |x| - 1 & \text{if $|x|\geq |y|$,} \\
    |y| - 1 & \text{else.}
  \end{cases}
\]
The phase-field function is then calculated directly from the signed-distance
function in \cref{eq:characteristic}.

\subsubsection{Case 3}

Again let $D$ be the circle centred at $(0,0)$ with radius $R=1$, but now
consider the case where the analytic solution is
\[
  u_{\text{an}}(x,y) = y\sqrt{x^2 + y^2},
\]
which corresponds to
\[
  f = \frac{3y}{\sqrt{x^2+y^2}} - y\sqrt{x^2 + y^2},
\]
$g = 2y$, and
\[
  \lapls u_{\text{an}} = -\frac{y}{\sqrt{x^2 + y^2}}.
\]
Note that in the DDM equations, $g$ is extrapolated constantly in the normal
direction off of the boundary $\partial D$.

\subsubsection{Case 4}

For the final Neumann case we again let $D=[-1,1]^2$, and we consider the case
where the analytic solution is
\[
  u_{\text{an}}(x,y) = e^r,
\]
where $r = \frac{x^2 + y^2}{4}$.  This corresponds to
\[
  f = r e^r.
\]
The boundary function $g$ and the surface Laplacian of the analytic function
along the boundary are
\begin{align*}
  g &= \frac 1 2 e^{\frac{1 + \xi^2}{4}}, \\
  \lapls u_\text{an} &= \frac 1 4\left(\xi^2 + 2\right)e^{\frac{1 + \xi^2}{4}},
\end{align*}
where $\xi \equiv x$ along the bottom and top boundaries, and $\xi \equiv y$
along the left and right boundaries.

\subsubsection{Results}

\Cref{fig:neumann_case1-plot,fig:neumann_case2-plot,fig:neumann_case3-plot,%
  fig:neumann_case4-plot} and \cref{tab:neumann_error_per_eps} show convergence
results in the $L^2$ norm where~$\epsilon$ is reduced for Cases 1 to 4 with
BC1, BC1M1, BC1M2 and BC2.
\Cref{fig:neumann_case4-plot-2,tab:neumann_error_per_eps_infty} shows the
results for Case 4 where the $L^\infty$ norm is used.  Although the DDM is most
efficient when adaptive meshes are used, here we consider only uniform meshes
to more easily control the discretization errors in order to focus on the
errors in the DDM.  As in all diffuse-interface methods, fine grids are
necessary to accurately solve the equations when $\epsilon$ is small.  This is
especially apparent for the cases with BC2 where even the finest grid spacing,
$n=8192$ in each direction, becomes too coarse to obtain results for small
$\epsilon$ that have converged with respect to the grid refinements.

The results confirm the second-order accuracy of all the considered boundary
condition approximations.  Note that while the difference between BC1 and BC1M1
tends to be small, BC1M1 consistently performs better than BC1.  In turn, BC1
performs better than BC2.  In case 2 there is a noticeable improvement of
BC1M1 over BC1.  Case 3 is the first case that has a nonconstant boundary
condition, and the surface Laplacian of the analytic solution along the
boundary is also nonconstant.  An unexpected result for Case 3 is that BC1M2
performs the best.  One possible explanation to this is errors due to grid
anisotropy.  Therefore we also consider a fourth case, which again has
a nonconstant boundary condition and nonconstant surface Laplacian of the
analytic solution.  Since the domain in this case is a square, the effect of
grid anisotropy is lessened.  Correspondingly BC1M1 performs the best.
The cases were also calculated with the $L^\infty$ norm, which gave similar
results, although at the smallest values of $\epsilon$, the orders of accuracy
of BC1M1 and BC1M2 may deteriorate in $L^\infty$, as seen in
\cref{fig:neumann_case4-plot-2,tab:neumann_error_per_eps_infty} for case
4.  This could be due to the influence of higher order terms in the expansion,
    or the amplification of error when $\epsilon$ is small due to the condition
    number of the system, which should scale like $\epsilon^{-2}$. This is
    currently under investigation.

The difference between BC1 and BC2 is noticeable, especially with regard to the
required amount of grid refinement that is needed to obtain a convergent
result.  This provides practical limits for the use of BC2.

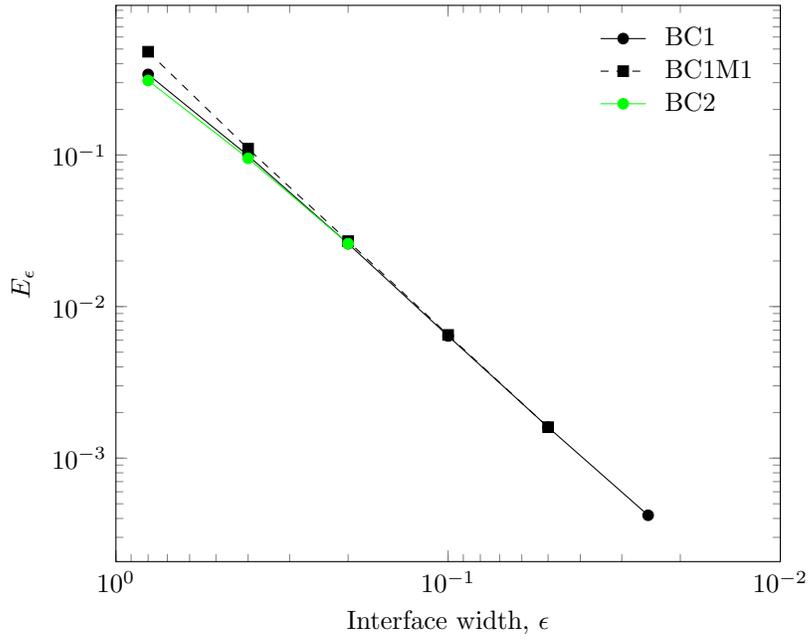
\begin{figure}[tbp]
  \centering
  \begin{tikzpicture}
    \begin{loglogaxis}[
      xlabel={Interface width, $\epsilon$},
      ylabel={$E_\epsilon$},
      x dir=reverse,
      width=0.8\textwidth,
      legend entries={BC1,
                      BC1M1,
                      BC2},
      legend cell align=left,
      legend style={column sep=0.5em,draw=white},
      ]

    \addplot[solid,mark=*,black] plot coordinates {
    (8.000e-01, 3.4e-01)
    (4.000e-01, 9.9e-02)
    (2.000e-01, 2.6e-02)
    (1.000e-01, 6.4e-03)
    (5.000e-02, 1.6e-03)
    (2.500e-02, 4.2e-04)
    };

    \addplot[dashed,mark=square*,mark options=solid,black] plot coordinates {
    (8.000e-01, 4.8e-01)
    (4.000e-01, 1.1e-01)
    (2.000e-01, 2.7e-02)
    (1.000e-01, 6.5e-03)
    (5.000e-02, 1.6e-03)
    };

    \addplot[solid,mark=*,green] plot coordinates {
    (8.000e-01, 3.1e-01)
    (4.000e-01, 9.5e-02)
    (2.000e-01, 2.6e-02)
    };

    \end{loglogaxis}
  \end{tikzpicture}
  \caption{$L^2$ errors for the Neumann problem with respect to $\epsilon$ for
    Case 1, as labelled.}
  \label{fig:neumann_case1-plot}
\end{figure}
\begin{figure}[tbp]
  \centering
  \begin{tikzpicture}
    \begin{loglogaxis}[
      xlabel={Interface width, $\epsilon$},
      ylabel={$E_\epsilon$},
      x dir=reverse,
      width=0.8\textwidth,
      legend entries={BC1,
                      BC1M1,
                      BC1M2,
                      BC2},
      legend cell align=left,
      legend style={column sep=0.5em,draw=white},
      ]

    \addplot[solid,mark=*,black] plot coordinates {
      (8.000e-01, 2.5e-01)
      (4.000e-01, 7.3e-02)
      (2.000e-01, 1.9e-02)
      (1.000e-01, 5.2e-03)
    };

    \addplot[dashed,mark=square*,mark options=solid,black] plot coordinates {
      (8.000e-01, 2.0e-01)
      (4.000e-01, 2.6e-02)
      (2.000e-01, 5.2e-03)
      (1.000e-01, 1.2e-03)
    };

    \addplot[dotted,mark=triangle*,mark options=solid,black] plot coordinates {
      (8.000e-01, 2.9e-01)
      (4.000e-01, 7.8e-02)
      (2.000e-01, 2.0e-02)
      (1.000e-01, 5.1e-03)
    };

    \addplot[solid,mark=*,green] plot coordinates {
      (8.000e-01, 2.2e-01)
      (4.000e-01, 7.0e-02)
      (2.000e-01, 2.0e-02)
    };

    \end{loglogaxis}
  \end{tikzpicture}
  \caption{$L^2$ errors for the Neumann problem with respect to $\epsilon$ for
  Case 2, as labelled.}
  \label{fig:neumann_case2-plot}
\end{figure}
\begin{figure}[tbp]
  \centering
  \begin{tikzpicture}
    \begin{loglogaxis}[
      xlabel={Interface width, $\epsilon$},
      ylabel={$E_\epsilon$},
      x dir=reverse,
      width=0.8\textwidth,
      legend entries={BC1,
                      BC1M1,
                      BC1M2,
                      BC2},
      legend cell align=left,
      legend style={column sep=0.5em,draw=white},
      ]

    \addplot[solid,mark=*,black] plot coordinates {
      (8.000e-01, 1.3e-01)
      (4.000e-01, 3.1e-02)
      (2.000e-01, 7.5e-03)
      (1.000e-01, 1.8e-03)
      (5.000e-02, 4.5e-04)
      (2.500e-02, 1.2e-04)
    };

    \addplot[dashed,mark=square*,mark options=solid,black] plot coordinates {
      (8.000e-01, 8.7e-02)
      (4.000e-01, 2.9e-02)
      (2.000e-01, 7.1e-03)
      (1.000e-01, 1.8e-03)
      (5.000e-02, 4.3e-04)
      (2.500e-02, 1.1e-04)
    };

    \addplot[dotted,mark=triangle*,mark options=solid,black] plot coordinates {
      (8.000e-01, 3.2e-02)
      (4.000e-01, 1.3e-02)
      (2.000e-01, 3.4e-03)
      (1.000e-01, 8.6e-04)
      (5.000e-02, 2.1e-04)
    };

    \addplot[solid,mark=*,green] plot coordinates {
      (8.000e-01, 1.2e-01)
      (4.000e-01, 2.8e-02)
      (2.000e-01, 8.1e-03)
    };

    \end{loglogaxis}
  \end{tikzpicture}
  \caption{$L^2$ errors for the Neumann problem with respect to $\epsilon$ for
    Case 3, as labelled.}
  \label{fig:neumann_case3-plot}
\end{figure}
\begin{figure}[tbp]
  \centering
  \begin{tikzpicture}
    \begin{loglogaxis}[
      xlabel={Interface width, $\epsilon$},
      ylabel={$E_\epsilon$},
      x dir=reverse,
      width=0.8\textwidth,
      legend entries={BC1,
                      BC1M1,
                      BC1M2,
                      BC2},
      legend cell align=left,
      legend style={column sep=0.5em,draw=white},
      ]

    \addplot[solid,mark=*,black] plot coordinates {
      (8.000e-01, 1.7e-01)
      (4.000e-01, 4.4e-02)
      (2.000e-01, 1.1e-02)
      (1.000e-01, 3.0e-03)
    };

    \addplot[dashed,mark=square*,mark options=solid,black] plot coordinates {
      (8.000e-01, 1.4e-01)
      (4.000e-01, 2.0e-02)
      (2.000e-01, 4.6e-03)
      (1.000e-01, 1.2e-03)
    };

    \addplot[dotted,mark=triangle*,mark options=solid,black] plot coordinates {
      (8.000e-01, 3.4e-01)
      (4.000e-01, 5.5e-02)
      (2.000e-01, 1.2e-02)
      (1.000e-01, 3.1e-03)
    };

    \addplot[solid,mark=*,green] plot coordinates {
      (8.000e-01, 1.7e-01)
      (4.000e-01, 4.6e-02)
      (2.000e-01, 1.2e-02)
    };

    \end{loglogaxis}
  \end{tikzpicture}
  \caption{$L^2$ errors for the Neumann problem with respect to $\epsilon$ for
    Case 4, as labelled.}
  \label{fig:neumann_case4-plot}
\end{figure}
\begin{figure}[tbp]
  \centering
  \begin{tikzpicture}
    \begin{loglogaxis}[
      xlabel={Interface width, $\epsilon$},
      ylabel={$E_\epsilon$},
      x dir=reverse,
      width=0.8\textwidth,
      legend entries={BC1,
                      BC1M1,
                      BC1M2,
                      BC2},
      legend cell align=left,
      legend style={column sep=0.5em,draw=white},
      ]

    \addplot[solid,mark=*,black] plot coordinates {
      (8.000e-01, 1.8e-01)
      (4.000e-01, 5.2e-02)
      (2.000e-01, 1.5e-02)
      (1.000e-01, 4.3e-03)
    };

    \addplot[dashed,mark=square*,mark options=solid,black] plot coordinates {
      (8.000e-01, 1.5e-01)
      (4.000e-01, 2.2e-02)
      (2.000e-01, 1.2e-02)
      (1.000e-01, 5.8e-03)
      (5.000e-02, 2.8e-03)
    };

    \addplot[dotted,mark=triangle*,mark options=solid,black] plot coordinates {
      (8.000e-01, 3.5e-01)
      (4.000e-01, 6.3e-02)
      (2.000e-01, 1.6e-02)
      (1.000e-01, 4.6e-03)
      (5.000e-02, 2.4e-03)
    };

    \addplot[solid,mark=*,green] plot coordinates {
      (8.000e-01, 1.8e-01)
      (4.000e-01, 5.1e-02)
      (2.000e-01, 1.4e-02)
      (1.000e-01, 4.7e-03)
    };

    \end{loglogaxis}
  \end{tikzpicture}
  \caption{$L^\infty$ errors for the Neumann problem with respect to $\epsilon$
    for Case 4, as labelled.}
  \label{fig:neumann_case4-plot-2}
\end{figure}
\begin{table}[tbp]
  \scriptsize
  \centering
  \begin{tabular}{rlrlrlrlr}
    \toprule
               & BC1 & &  BC1M1 & &  BC1M2 & & BC2 \\
    $\epsilon$ & $E$ & $k$ & $E$ & $k$ & $E$ & $k$ & $E$ & $k$ \\
    \midrule
    \multicolumn{4}{l}{Case 1} \\
    0.800 & 3.39\e 1&     & 4.77\e 1&     &         &     & 3.09\e 1&     \\
    0.400 & 9.94\e 2& 1.8 & 1.12\e 1& 2.1 &         &     & 9.52\e 2& 1.7 \\
    0.200 & 2.57\e 2& 2.0 & 2.68\e 2& 2.1 &         &     & 2.57\e 2& 1.9 \\
    0.100 & 6.43\e 3& 2.0 & 6.51\e 3& 2.0 &         &     &         &     \\
    0.050 & 1.61\e 3& 2.0 & 1.59\e 3& 2.0 &         &     &         &     \\
    0.025 & 4.15\e 4& 2.0 & 3.87\e 4& 2.0 &         &     &         &     \\
    \midrule
    \multicolumn{4}{l}{Case 2} \\
    0.800 & 2.46\e 1&     & 1.96\e 1&     & 2.88\e 1&     & 2.22\e 1&     \\
    0.400 & 7.30\e 2& 1.8 & 2.58\e 2& 2.9 & 7.81\e 2& 1.9 & 6.99\e 2& 1.7 \\
    0.200 & 1.94\e 2& 1.9 & 5.21\e 3& 2.3 & 1.96\e 2& 2.0 & 1.95\e 2& 1.8 \\
    0.100 & 5.16\e 3& 1.9 & 1.20\e 3& 2.1 & 5.10\e 3& 1.9 &         &     \\
    \midrule
    \multicolumn{4}{l}{Case 3} \\
    0.800 & 1.27\e 1&     & 8.74\e 2&     & 3.16\e 2&     & 1.18\e 1&     \\
    0.400 & 3.12\e 2& 2.0 & 2.85\e 2& 1.6 & 1.28\e 2& 1.3 & 2.82\e 2& 2.1 \\
    0.200 & 7.48\e 3& 2.1 & 7.08\e 3& 2.0 & 3.40\e 3& 1.9 & 8.13\e 3& 1.8 \\
    0.100 & 1.81\e 3& 2.0 & 1.75\e 3& 2.0 & 8.58\e 4& 2.0 &         &     \\
    0.050 & 4.48\e 4& 2.0 & 4.32\e 4& 2.0 & 2.12\e 4& 2.0 &         &     \\
    0.025 & 1.15\e 4& 2.0 & 1.06\e 4& 2.0 &         &     &         &     \\
    \midrule
    \multicolumn{4}{l}{Case 4} \\
    0.800 & 1.71\e 1&     & 1.38\e 1&     & 3.39\e 1&     & 1.74\e 1&     \\
    0.400 & 4.42\e 2& 2.0 & 2.04\e 2& 2.8 & 5.51\e 2& 2.6 & 4.61\e 2& 1.9 \\
    0.200 & 1.14\e 2& 2.0 & 4.58\e 3& 2.2 & 1.24\e 2& 2.2 & 1.20\e 2& 1.9 \\
    0.100 & 2.95\e 3& 1.9 & 1.18\e 3& 2.0 & 3.09\e 3& 2.0 &         &     \\
    \bottomrule
  \end{tabular}
  \caption{The $L^2$ error for the Neumann problem as a function of $\epsilon$
    for all cases.  All results are calculated with $n=8192$ in each direction
    on uniform grids.  Except for the scheme BC1M2 in Case 1, which was not
    simulated, blank results indicate that the solutions require even finer
    grids to converge.}
  \label{tab:neumann_error_per_eps}
\end{table}
\begin{table}[tbp]
  \scriptsize
  \centering
  \begin{tabular}{rlrlrlrlr}
    \toprule
               & BC1 &     & BC1M1 &     & BC1M2 &     & BC2 &     \\
    $\epsilon$ & $E$ & $k$ & $E$   & $k$ & $E$   & $k$ & $E$ & $k$ \\
    \midrule
    0.800 & 1.80\e 1&     & 1.46\e 1&     & 3.54\e 1&     & 1.76\e 1&     \\
    0.400 & 5.17\e 2& 1.8 & 2.24\e 2& 2.7 & 6.32\e 2& 2.5 & 5.08\e 2& 1.8 \\
    0.200 & 1.48\e 2& 1.8 & 1.18\e 2& 0.9 & 1.59\e 2& 2.0 & 1.44\e 2& 1.8 \\
    0.100 & 4.29\e 3& 1.8 & 5.77\e 3& 1.0 & 4.56\e 3& 1.8 & 4.72\e 3& 1.6 \\
    0.050 &         &     & 2.79\e 3& 1.0 & 2.40\e 3& 0.9 &         &     \\
    \bottomrule
  \end{tabular}
  \caption{The $L^\infty$ error for the Neumann problem as a function of
    $\epsilon$ for Case 4.  All results are calculated with $n=8192$ in each
    direction on uniform grids.  Blank results indicate that the solutions
    require even finer grids to converge.}
  \label{tab:neumann_error_per_eps_infty}
\end{table}

\clearpage
\subsection{Robin boundary conditions}

Now consider the Poisson equation with Robin boundary conditions,
\[
  \begin{alignedat}{3}
    \lapl u      &= f            && \txin D, \\
    \ndot\grad u &= k(u-g) \quad && \txon \partial D.
  \end{alignedat}
\]
As in the previous section, we solve the DDM equation
\[
  \div\left(\phi\grad u\right) + \text{BC} = \phi f,
\]
using BC1, BC2, BC1M and BC2M.

\subsubsection{Case 1}

Consider the case where $D$ is a circle of radius $R=1$ centred at $(0,0)$,
and where the analytic solution to the Poisson equation in $D$ is
\[
  u_{\text{an}}(x,y) = \frac{1}{4}\left( x^2 + y^2 \right).
\]
This corresponds to $f = 1$,
\[
  g = \frac{1}{2} \left(\frac{1}{2} - \frac{1}{k}\right),
\]
and $\lapls u_{\text{an}} = 0$.  We will consider the case when $k=-1$, thus
$g=3/4$.

\subsubsection{Case 2}

Again let $D$ be the circle at $(0,0)$ with radius $R=1$, but now consider the
case where the analytic solution is
\[
  u_{\text{an}}(x,y) = y\left(x^2 + y^2\right),
\]
which corresponds to
\begin{align*}
  f &= -8y, \\
  g &= y\left( 1 - \frac 3 k \right),
\end{align*}
and
\[
  \lapls u_{\text{an}} = -y,
\]
Again let $k=-1$ so that $g=4y$.  Similar to the Neumann case 3, $g$ is
extended constantly in the normal direction in the DDM equations.

\subsubsection{Case 3}

For the final Robin case we let $D=[-1,1]^2$, and we consider a case that
corresponds to the Neumann Case 4 where the analytic solution is
\[
  u_{\text{an}}(x,y) = e^r,
\]
where $r = \frac{x^2 + y^2}{4}$.  This corresponds to
\[
  f = (r+1) e^r.
\]
The boundary function $g$ and the surface Laplacian of the analytic function
along the boundary are
\begin{align*}
  g &= \frac 3 2 e^{\frac{1 + \xi^2}{4}}, \\
  \lapls u_\text{an} &= \frac 1 4\left(\xi^2 + 2\right)e^{\frac{1 + \xi^2}{4}},
\end{align*}
where $\xi \equiv x$ along the bottom and top boundaries, and $\xi \equiv y$
along the left and right boundaries.

\subsubsection{Results}

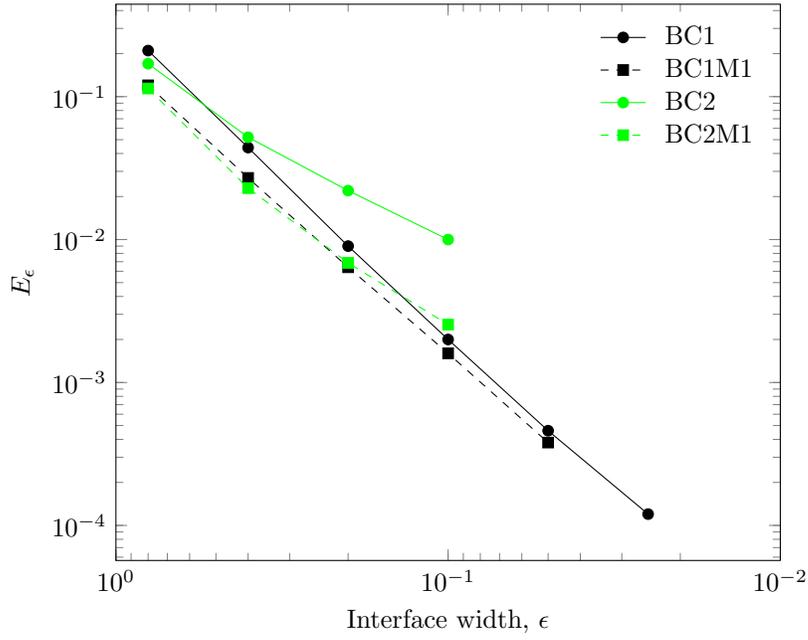
\begin{figure}[tbp]
  \centering
  \begin{tikzpicture}
    \begin{loglogaxis}[
      xlabel={Interface width, $\epsilon$},
      ylabel={$E_\epsilon$},
      x dir=reverse,
      width=0.8\textwidth,
      legend entries={BC1,
                      BC1M1,
                      BC2,
                      BC2M1,},
      legend cell align=left,
      legend style={column sep=0.5em,draw=white},
      ]

    \addplot[solid,mark=*,black] plot coordinates {
      (8.000e-01, 2.1e-01)
      (4.000e-01, 4.4e-02)
      (2.000e-01, 9.0e-03)
      (1.000e-01, 2.0e-03)
      (5.000e-02, 4.6e-04)
      (2.500e-02, 1.2e-04)
    };

    \addplot[dashed,mark=square*,mark options=solid,black] plot coordinates {
      (8.000e-01, 1.2e-01)
      (4.000e-01, 2.7e-02)
      (2.000e-01, 6.4e-03)
      (1.000e-01, 1.6e-03)
      (5.000e-02, 3.8e-04)
    };

    \addplot[solid,mark=*,green] plot coordinates {
      (8.000e-01, 1.7e-01)
      (4.000e-01, 5.2e-02)
      (2.000e-01, 2.2e-02)
      (1.000e-01, 1.0e-02)
    };

    \addplot[dashed,mark=square*,mark options=solid,green] plot coordinates {
      (8.000e-01, 1.14e-01)
      (4.000e-01, 2.29e-02)
      (2.000e-01, 6.87e-03)
      (1.000e-01, 2.54e-03)
    };

    \end{loglogaxis}
  \end{tikzpicture}
  \caption{$L^2$ errors for the Robin problem with respect to $\epsilon$ for
    Case 1, as labelled.}
  \label{fig:robin_cases-plot1}
\end{figure}

The convergence results calculated with the $L^2$ norm are presented in
\cref{fig:robin_cases-plot1,fig:robin_cases-plot2,fig:robin_cases-plot3} and
\cref{tab:robin_error_per_eps}.
\Cref{fig:robin_cases-plot2-2,tab:robin_error_per_eps_infty} shows the results
for Case 2 where the $L^\infty$ norm is used.  Again the results indicate that
BC1M1 performs better than BC1, although both methods are second-order
accurate, as predicted by our analysis.  The results also show that BC1 gives
better results than BC2, which is approximately first-order accurate as also
predicted by theory.  Further, as in the Neumann case, BC2 is seen to require
very fine grids to converge.  For small $\epsilon$, the requirement exceeds our
finest grid.

The modified BC2M1 and BC2M2 schemes are also tested.  The results with BC2M2
are almost indistinguishable from the results with BC2M1, so only the latter
results are shown in the following figures.  All results are listed in
\cref{tab:robin_error_per_eps} and \cref{tab:robin_error_per_eps_infty}.  The
BC2M schemes are shown to perform better than the BC2 scheme, but they also
require very fine grids to converge.  Further, the orders of accuracy of BC2M1
and BC2M2 seem to deteriorate somewhat at the smallest values of $\epsilon$.
As discussed in \cref{sec:dda_robin_bc2m}, this could be due to the influence
of higher order terms in the expansion, or the amplification of error and is
under study.

\begin{figure}[tbp]
  \centering
  \begin{tikzpicture}
    \begin{loglogaxis}[
      xlabel={Interface width, $\epsilon$},
      ylabel={$E_\epsilon$},
      x dir=reverse,
      width=0.8\textwidth,
      legend entries={BC1,
                      BC1M1,
                      BC2,
                      BC2M1,},
      legend cell align=left,
      legend style={column sep=0.5em,draw=white},
      ]

    \addplot[solid,mark=*,black] plot coordinates {
      (8.000e-01, 4.21e-01)
      (4.000e-01, 1.02e-01)
      (2.000e-01, 2.19e-02)
      (1.000e-01, 4.77e-03)
      (5.000e-02, 1.09e-03)
      (2.500e-02, 2.57e-04)
    };

    \addplot[dashed,mark=square*,mark options=solid,black] plot coordinates {
      (8.000e-01, 1.27e-1)
      (4.000e-01, 4.87e-2)
      (2.000e-01, 1.30e-2)
      (1.000e-01, 3.27e-3)
      (5.000e-02, 8.15e-4)
      (2.500e-02, 2.04e-4)
    };

    \addplot[solid,mark=*,green] plot coordinates {
      (8.000e-01, 3.84e-1)
      (4.000e-01, 9.16e-2)
      (2.000e-01, 2.51e-2)
      (1.000e-01, 9.13e-3)
    };

    \addplot[dashed,mark=square*,mark options=solid,green] plot coordinates {
      (8.000e-01, 3.48e-1)
      (4.000e-01, 7.24e-2)
      (2.000e-01, 1.56e-2)
      (1.000e-01, 4.65e-3)
    };

    \end{loglogaxis}
  \end{tikzpicture}
  \caption{$L^2$ errors for the Robin problem with respect to $\epsilon$ for
    Case 2, as labelled.}
  \label{fig:robin_cases-plot2}
\end{figure}
\begin{figure}[tbp]
  \centering
  \begin{tikzpicture}
    \begin{loglogaxis}[
      xlabel={Interface width, $\epsilon$},
      ylabel={$E_\epsilon$},
      x dir=reverse,
      width=0.8\textwidth,
      legend entries={BC1,
                      BC1M1,
                      BC2,
                      BC2M1,},
      legend cell align=left,
      legend style={column sep=0.5em,draw=white},
      ]

    \addplot[solid,mark=*,black] plot coordinates {
      (8.000e-01, 3.80e-1)
      (4.000e-01, 9.39e-2)
      (2.000e-01, 2.14e-2)
      (1.000e-01, 4.90e-3)
      (0.500e-01, 1.15e-3)
      (0.250e-01, 2.78e-4)
    };

    \addplot[dashed,mark=square*,mark options=solid,black] plot coordinates {
      (8.000e-01, 1.36e-1)
      (4.000e-01, 4.30e-2)
      (2.000e-01, 1.06e-2)
      (1.000e-01, 2.53e-3)
      (5.000e-02, 6.09e-4)
      (2.500e-02, 1.49e-4)
    };

    \addplot[solid,mark=*,green] plot coordinates {
      (8.000e-01, 3.24e-1)
      (4.000e-01, 7.49e-2)
      (2.000e-01, 1.91e-2)
      (1.000e-01, 6.98e-3)
    };

    \addplot[dashed,mark=square*,mark options=solid,green] plot coordinates {
      (8.000e-01, 3.04e-1)
      (4.000e-01, 6.94e-2)
      (2.000e-01, 1.72e-2)
      (1.000e-01, 6.40e-3)
    };

    \end{loglogaxis}
  \end{tikzpicture}
  \caption{$L^\infty$ errors for the Robin problem with respect to $\epsilon$
    for Case 2, as labelled.}
  \label{fig:robin_cases-plot2-2}
\end{figure}
\begin{figure}[tbp]
  \centering
  \begin{tikzpicture}
    \begin{loglogaxis}[
      xlabel={Interface width, $\epsilon$},
      ylabel={$E_\epsilon$},
      x dir=reverse,
      width=0.8\textwidth,
      legend entries={BC1,
                      BC1M1,
                      BC2,
                      BC2M1,},
      legend cell align=left,
      legend style={column sep=0.5em,draw=white},
      ]

    \addplot[solid,mark=*,black] plot coordinates {
      (8.000e-01, 7.9e-02)
      (4.000e-01, 1.6e-02)
      (2.000e-01, 3.7e-03)
      (1.000e-01, 9.0e-04)
    };

    \addplot[dashed,mark=square*,mark options=solid,black] plot coordinates {
      (8.000e-01, 3.6e-02)
      (4.000e-01, 7.4e-03)
      (2.000e-01, 1.7e-03)
      (1.000e-01, 4.3e-04)
    };

    \addplot[solid,mark=*,green] plot coordinates {
      (8.000e-01, 8.2e-01)
      (4.000e-01, 1.9e-02)
      (2.000e-01, 5.9e-03)
    };

    \addplot[dashed,mark=square*,mark options=solid,green] plot coordinates {
      (8.000e-01, 6.75e-02)
      (4.000e-01, 1.27e-02)
      (2.000e-01, 2.78e-03)
    };

    \end{loglogaxis}
  \end{tikzpicture}
  \caption{$L^2$ errors for the Robin problem with respect to $\epsilon$ for
    Case 3, as labelled.}
  \label{fig:robin_cases-plot3}
\end{figure}
\begin{table}[tbp]
  \tiny
  \centering
  \begin{tabular}{clrlrlrlrlr}
    \toprule
               & BC1 & & BC1M1 & & BC2 & & BC2M1 & & BC2M2\\
    $\epsilon$ & $E$ & $k$ & $E$ & $k$ & $E$ & $k$ & $E$ & $k$ & $E$ & $k$ \\
    \midrule
    \multicolumn{4}{l}{Case 1} \\
    0.800&2.11\e 1&    &1.20\e 1&    &1.70\e 1&   &1.14\e 1&   &1.15\e 1\\
    0.400&4.40\e 2&2.3 &2.72\e 2&2.1 &5.21\e 2&1.7&2.29\e 2&2.3&2.31\e 2 &2.3\\
    0.200&8.99\e 3&2.3 &6.42\e 3&2.1 &2.18\e 2&1.3&6.87\e 3&1.7&6.92\e 3 &1.7\\
    0.100&1.95\e 3&2.2 &1.57\e 3&2.0 &1.03\e 2&1.1&2.54\e 3&1.4&2.55\e 3 &1.4\\
    0.050&4.57\e 4&2.1 &3.79\e 4&2.0\\
    0.025&1.23\e 4&1.9\\
    \midrule
    \multicolumn{4}{l}{Case 2} \\
    0.800&4.21\e 1&   &1.27\e 1&   &3.84\e 1&    &3.48\e 1&   &3.49\e 1&   \\
    0.400&1.02\e 1&2.0&4.87\e 2&1.4&9.16\e 2&2.1 &7.24\e 2&2.3&7.26\e 2&2.3\\
    0.200&2.19\e 2&2.2&1.30\e 2&1.9&2.51\e 2&1.9 &1.56\e 2&2.2&1.57\e 2&2.2\\
    0.100&4.77\e 3&2.2&3.27\e 3&2.0&9.13\e 3&1.5 &4.65\e 3&1.7&4.67\e 3&1.7\\
    0.050&1.09\e 3&2.1&8.15\e 4&2.0\\
    0.025&2.57\e 4&2.1&2.04\e 4&2.0\\
    \midrule
    \multicolumn{4}{l}{Case 3} \\
    0.800&7.89\e 2&   &3.60\e 2&    &8.23\e 2&   &6.75\e 2&   &6.80\e 2&   \\
    0.400&1.64\e 2&2.3&7.38\e 3&2.3 &1.89\e 2&2.2&1.27\e 2&2.3&1.28\e 2&2.4\\
    0.200&3.70\e 3&2.2&1.71\e 3&2.1 &5.90\e 3&1.7&2.78\e 3&2.2&2.81\e 3&2.2\\
    0.100&9.04\e 4&2.0&4.28\e 4&2.0\\
    \bottomrule
  \end{tabular}
  \caption{The $L^2$ error for the Robin problem as a function of $\epsilon$
    for all cases.  All results are calculated with $n=8192$ in each direction
    on uniform grids, except for Case 3 with BC2, where the results are
    calculated with $n=4096$ in each direction.  Blank results indicate that
    the solutions require even finer grids to converge.}
  \label{tab:robin_error_per_eps}
\end{table}
\begin{table}[tbp]
  \tiny
  \centering
  \begin{tabular}{clrlrlrlrlr}
    \toprule
    & BC1 & & BC1M1 & & BC2 & & BC2M1 & & BC2M2 \\
    $\epsilon$ & $E$ & $k$ & $E$ & $k$ & $E$ & $k$ & $E$ & $k$ & $E$ & $k$ \\
    \midrule
    \multicolumn{4}{l}{Case 2} \\
    0.800&3.80\e 1&   &1.36\e 1&   &3.24\e 1&    &3.04\e 1&   &3.05\e 1&   \\
    0.400&9.39\e 2&2.0&4.30\e 2&1.7&7.49\e 2&2.1 &6.94\e 2&2.1&6.95\e 2&2.1\\
    0.200&2.14\e 2&2.1&1.06\e 2&2.0&1.91\e 2&2.0 &1.72\e 2&2.0&1.73\e 2&2.0\\
    0.100&4.90\e 3&2.1&2.53\e 3&2.1&6.98\e 3&1.5 &6.40\e 3&1.4&6.42\e 3&1.4\\
    0.050&1.15\e 3&2.1&6.09\e 4&2.1\\
    0.025&2.78\e 4&2.1&1.49\e 4&2.0\\
    \bottomrule
  \end{tabular}
  \caption{The $L^\infty$ error for the Robin problem as a function of
    $\epsilon$ for Case 2.  All results are calculated with $n=8192$ in each
    direction on uniform grids.  Blank results indicate that the solutions
    require even finer grids to converge.}
  \label{tab:robin_error_per_eps_infty}
\end{table}

\begin{figure}[tbp]
  \centering
  \begin{subfigure}{\textwidth}
    \centering
    \includegraphics[width=0.8\textwidth]{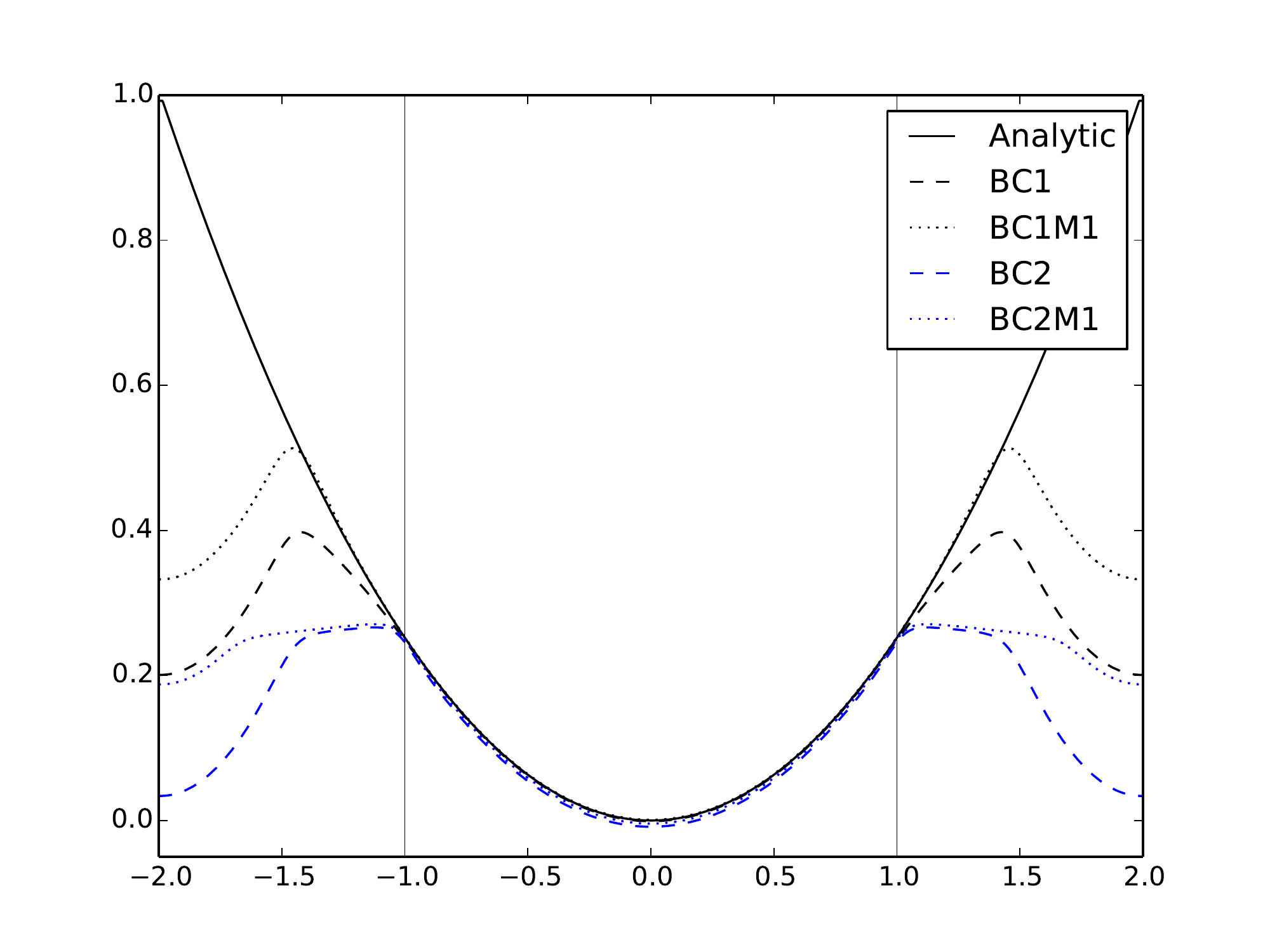}
    \caption{}
    \label{fig:plot_full}
  \end{subfigure}
  \begin{subfigure}{\textwidth}
    \centering
    \includegraphics[width=0.8\textwidth]{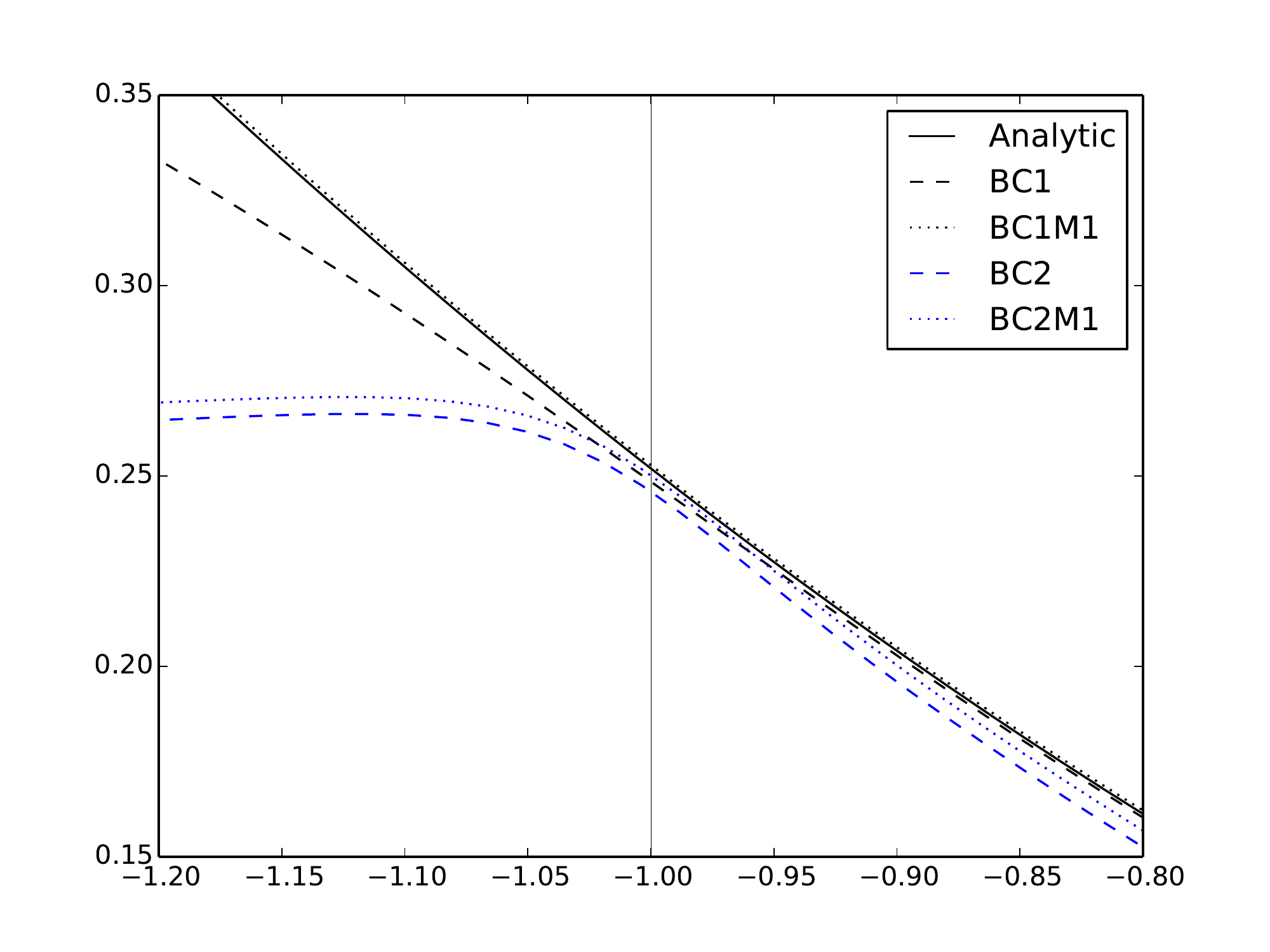}
    \caption{}
    \label{fig:plot_zoom}
  \end{subfigure}
  \caption{A plot of the solutions of Case 1 through $y=0$.  The solutions
    with BC2M1 and BC2M2 are indistinguishable, so only BC2M1 is shown.  (a)
    Shows the full slice, where the domain boundary is depicted as thin
    vertical lines at $x=\pm1$.  (b) A zoom-in that shows the solutions near
    the left boundary.}
  \label{fig:plot_robin_case1}
\end{figure}

\Cref{fig:plot_robin_case1} shows a plot of the solutions of Case 1 with
$\epsilon=0.2$ at $y=0$.  The plot shows the solutions with BC1 (black dashed),
BC1M1 (black doted), BC2 (blue dashed) and BC2M1 (blue dotted).  We see that
the solutions with the modified schemes BC1M1 and BC2M1 perform better than the
corresponding schemes with BC1 and BC2.

\clearpage
\section{Conclusion}
\label{sec:conclusion}

We have performed a matched asymptotic analysis of the DDM for the
Poisson equation with Robin boundary conditions and for a steady
reaction-diffusion equation with Neumann boundary conditions.  Our analysis
shows that for certain choices of the boundary condition approximations, the
DDM is second-order accurate in the interface thickness $\epsilon$.  However,
for other choices the DDM is only first-order accurate.  This is confirmed
numerically and helps to explain why the choice of boundary-condition
approximation is important for rapid global convergence and high accuracy.
This helps to explain why the choice of boundary-condition approximation is
important for rapid global convergence and high accuracy.  In particular, the
boundary condition BC1, which arises from representing the surface delta
function as $|\nabla\phi|$, is seen to give rise to a second-order
approximation for both the Neumann and Robin boundary conditions and thus is
perhaps the most reliable choice.  The boundary condition BC2, which arises
from approximating the surface delta function as $\epsilon|\nabla\phi|^2$
yields a second-order accurate approximation for the Neumann problem, but
only first-order accuracy for the Robin problem.  In addition, BC2 requires
very fine meshes to converge.

Our analysis also suggests correction terms that may be added to yield a more
accurate diffuse-domain method.  We have presented several techniques for
obtaining second-order boundary conditions and performed numerical
simulations that confirm the predicted accuracy, although the order of
accuracy may deteriorate at the smallest values of $\epsilon$ possibly due to
amplification errors associated with conditioning of the system or the
influence of higher order terms in the asymptotic expansion.  This is
currently under study.  Further, the correction terms do not improve the mesh
requirements for convergence.

A common feature of the correction terms is that the interface thickness must
be sufficiently small in order for the DDM to remain an elliptic equation.
In addition, one choice of boundary condition involves the use of the surface
Laplacian of the solution, which could in principle lead to faster asymptotic
convergence since it directly cancels terms in the inner expansion of the
asymptotic matching.  However, the extension of this term outside the domain
of interest can cause the loss of ellipticity of the DDM.  As such, this is
an intriguing but not a practical scheme.  Nevertheless, as a proof of
principle, we still considered the effect of this term, however, by using the
surface Laplacian of the analytic solution in the DDM.  We found that this
choice gave the smallest errors in nearly all the cases considered.  By
incorporating different extensions of the boundary conditions in the exterior
of the domain that automatically guarantee ellipticity, we aim to make this
method practical.  This is the subject of future investigations.

We plan to extend our analysis to the Dirichlet problem where the boundary
condition approximations considered by Li et al.\ \cite{Li09} seem only to
yield first-order accuracy \cite{Franz12,Reuter12}.  Our asymptotic analysis
thus has the potential to identify correction terms that can be used to
generate second-order accurate diffuse-domain methods for the Dirichlet
problem.

\medskip
{\bf Acknowledgement.}
KYL acknowledges support from the Fulbright foundation for a Visiting
Researcher Grant to fund a stay at the University of California, Irvine.  KYL
also acknowledges support from Statoil and GDF SUEZ, and the Research Council
of Norway (193062/S60) for the research project Enabling low emission LNG
systems.  JL acknowledges support from the National Science Foundation,
Division of Mathematical Sciences, and the National Institute of Health through
grant P50GM76516 for a Center of Excellence in Systems Biology at the
University of California, Irvine.  The authors gratefully thank Bernhard Müller
(NTNU) and Svend Tollak Munkejord (SINTEF Energy Research) for helpful
discussions and for feedback on the manuscript.  The authors also wish to thank
the anonymous reviewers for comments that greatly improved the manuscript.

\medskip
\bibliographystyle{siam}

\begin{thebibliography}{10}

\bibitem{Aland10}
{\sc S.~Aland, J.~Lowengrub, and A.~Voigt}, {\em Two-phase flow in complex
  geometries: A diffuse-domain approach.}, Computer Modeling in Engineering \&
Sciences, 57 (2010), pp.~77--106.

\bibitem{Almgren99}
{\sc R.~F. Almgren}, {\em Second-order phase field asymptotics for unequal
  conductivities}, SIAM Journal on Applied Mathematics, 59 (1999),
pp.~2086--2107.

\bibitem{Bedrossian10}
{\sc J. Bedrossian, J. H. von Brecht, S. Zhu, E. Sifakis, and J. M. Teran},
{\em A second-order virtual node method for elliptic problems with interfaces
  and irregular domains}, J. Comput. Phys., 229 (2010), pp.~6405-6426.

\bibitem{bernauer12}
{\sc M.K. Bernauer, R. Herzog}, {\em Implementation of an X-FEM solver for the
  classical two-phase {S}tefan problem}, J. Sci. Comput., 52 (2012),
pp.~271-293.

\bibitem{Bueno06b}
{\sc A.~Bueno-Orovio and V.~M. Perez-Garcia}, {\em Spectral smoothed boundary
  methods: the role of external boundary conditions}, Numer. Meth. Partial
Diff. Eqns., 22 (2006), pp.~435--448.

\bibitem{Bueno06a}
{\sc A.~Bueno-Orovio, V.~M. Perez-Garcia, and F.~H. Fenton}, {\em Spectral
  methods for partial differential equations in irregular domains: the spectral
  smoothed boundary method}, SIAM Journal on Scientific Computing, 28 (2006),
pp.~886--900.

\bibitem{byfut12}
{\sc A. Byfut, A. Schroeder}, {\em hp-adaptive extended finite element method},
Int. J. Numer. Meth. Eng., 89 (2012), pp.~1293-1418.

\bibitem{Caginalp88}
{\sc G. Caginalp, P. C. Fife}, {\em Dynamics of Layered Interfaces Arising from
  Phase Boundaries}, SIAM Journal on Applied Mathematics, 48 (1988),
pp.~506--518.

\bibitem{cisternino12}
{\sc M. Cisternino, L Weynans}, {\em A parallel second-order Cartesian method
  for elliptic interface problems}, Comm. Comput. Phys., 12 (2012),
pp.~1562-1587.

\bibitem{Coco13}
{\sc A.~Coco and G.~Russo}, {\em Finite-Difference Ghost-Point Multigrid
  Methods on Cartesian Grids for Elliptic Problems in Arbitrary Domains},
J. Comput. Phys., 241 (2013), pp.~464-501.

\bibitem{DD07}
{\sc A. Demlow and G. Dziuk}, {\em An adaptive finite element method for the
  Laplace-Beltrami operator on implicitly defined surfaces}, SIAM J. Numer.
Anal., 45 (2007), pp.~421-442.

\bibitem{dolbow09}
{\sc J. Dolbow, I. Harari}, {\em An efficient finite element method for
  embedded interface problems}, Int. J. Numer. Meth. Eng., 78 (2009),
pp.~229-252.

\bibitem{duddu11}
{\sc R. Duddu, D.L. Chopp, P. Voorhees, B. Moran}, {\em Diffusional evolution
  of precipitates in elastic media using the extended finite element method and
  level set methods}, J. Comput. Phys., 230 (2011), pp.~1249-1264.

\bibitem{dziuk08a}
{\sc G. Dziuk and C.M. Elliott}, {\em Eulerian finite element method for
  parabolic PDEs on implicit surfaces}, Int. Free. Bound., 10 (2008),
pp.~119-138.

\bibitem{dziuk08b}
{\sc G. Dziuk and C.M. Elliott}, {\em An {E}ulerian approach to transport and
  diffusion on evolving implicit surfaces}, Comput. Visual. Sci., 13, (2010),
pp.~17-28.

\bibitem{DE12}{\sc G. Dziuk and C.M. Elliott}, {\em A fully discrete evolving
  surface finite element method}, SIAM J. Numer. Anal., 50, 5, (2012),
pp.~2677-2694.

\bibitem{elliott09}
{\sc C.M. Elliott,  B. Stinner,  V. Styles, R. Welford}, {\em Numerical
  computation of advection and diffusion on evolving diffuse interfaces}, IMA
J. Num. Anal., 31, (2011), pp.~245-269.

\bibitem{elliott09a}
{\sc C.M. Elliott and B. Stinner}, {\em Analysis of a diffuse interface
  approach to an advection diffusion equation on a moving surface}, Math. Mod.
Meth. Appl. Sci., (2009) in press.

\bibitem{FedkiwAslamMerrimanOsher_1999}
{\sc R.P. Fedkiw, T. Aslam, B. Merriman, S. Osher}, {\em A non-oscillatory
  {E}ulerian approach to interfaces in multimaterial flows (the ghost fluid
  method}, J. Comput. Phys., 152 (1999), pp.~457-492.

\bibitem{Fenton05}
{\sc F.~H. Fenton, E.~M. Cherry, A.~Karma, and W.~J. Rappel}, {\em Modeling
  wave propagation in realistic heart geometries using the phase-field method},
{CHAOS}, 15 (2005).

\bibitem{Folch1999}
{\sc R. Folch, J. Casademunt, A. Hernandez-Machado, and L. Ramirez-Piscina},
Phys. Rev. E, 60 (1999), pp.~1724.

\bibitem{Franz12}
{\sc S.~Franz, R.~Gärtner, H.-G.~Roos, and A.~Voigt}, {\em A Note on the
  Convergence Analysis of a Diffuse-Domain Approach},
Computational Methods in Applied Mathematics, 12 (2012), pp.~153--167.

\bibitem{fries10}
{\sc F.-P. Fries, T. Belytschko}, {\em The extended/generalized finite element
  method: An overview of the method and its applications}, Int. J. Numer. Meth.
Eng., 84 (2010), pp.~253-304.

\bibitem{Gibou13}
{\sc F. Gibou, C. Min, and R. Fedkiw}, {\em High Resolution Sharp Computational
  Methods for Elliptic and Parabolic Problems in Complex Geometries}, Journal
of Scientific Computing, 54 (2013), pp.~369-413.

\bibitem{GibouFedkiwChengKang_2002}
{\sc F. Gibou,  R. Fedkiw,  L.T. Cheng, and M. Kang}, {\em A second-order
  accurate symetric discretization of the {P}oisson equation on irregular
  domains}, J. Comput. Phys., 176 (2002), pp.~205-227.

\bibitem{GibouFedkiw_2005}
{\sc F. Gibou and R. Fedkiw}, {\em A fourth order accurate discretization for
  the {L}aplace and heat equations on arbitrary domains with applications to
  the {S}tefan problem}, J. Comput. Phys., 202 (2005), pp.~577-601.

\bibitem{GlimmMarchesinMcBryan_1981}
{\sc J. Glimm and D. Marchesin and O. McBryan}, {\em A numerical method for
  2 phase flow with an unstable interface}, J. Comput. Phys., 39 (1981),
pp.~179-200.

\bibitem{GlowinskiPanPeriaux_CMAME_1994}
{\sc R. Glowinski, T.W. Pan, and J. Periaux}, {\em A fictitious domain method
  for external incompressible viscous-flow modeled by {N}avier-{S}tokes
  equations}, Comput. Meth. Appl. Mech. Engin., 112 (1994), pp.~133-148.

\bibitem{GlowinskiPanWellsZhou_JCP_1996}
{\sc  R. Glowinski and T.W. Pan and R.O. Wells and X.D. Zhou}, {\em Wavelet and
  finite element solutions for the {N}eumann problem using fictitious domains},
J. Comput. Phys., 126 (1996), pp.~40-51.

\bibitem{greer06}
{\sc J.B. Greer and A.L. Bertozzi and G. Sapiro}, {\em Fourth order partial
  differential equations on general geometries}, J. Comput. Phys., 216 (2006),
pp.~216-246.

\bibitem{GR07}
{\sc S. Gross, and A. Reusken}, {\em An extended pressure finite element space
  for two-phase incompressible flows}, J. Comput. Phys., 224 (2007), pp.~40-48.

\bibitem{he11}
{\sc X.M. He, T. Lin, Y.P. Lin}, {\em Immersed finite element methods for
  elliptic interface problems with non-homogeneous jump conditions}, Int. J.
Numer. Anal. Model., 8 (2011), pp.~284-301.

\bibitem{Helgadottir11}
{\sc Á. Helgadóttir and F. Gibou}, {\em A Poisson–Boltzmann solver on irregular
  domains with Neumann or Robin boundary conditions on non-graded adaptive
  grid}, J. Comput. Phys., 230 (2011), pp.~3830-3848.

\bibitem{Hellrung12}
{\sc J. L. Hellrung Jr., L. Wang, E. Sifakis, and J. M. Teran}, {\em A second
  order virtual node method for elliptic problems with interfaces and irregular
  domains in three dimensions}, J. Comput. Phys., 231 (2012), pp.~2015-2048.

\bibitem{JiLienYee_2006}
{\sc H. Ji and F.-S. Lien and E. Yee}, {\em An efficient second-order accurate
  cut-cell method for solving the variable coefficient {P}oisson equation with
  jump conditions on irregular domains}, Int. J. Numer. Meth. Fluids, 52
(2006), pp.~723-748.

\bibitem{JohansenColella_1998}
{\sc H. Johansen and P. Colella}, {\em A {C}artesian grid embedded boundary
  method for {P}oisson's equation on irregular domains}, J. Comput. Phys., 147
(1998), pp.~60-85.

\bibitem{Colellaetal_2008}
{\sc H. Johansen and P. Colella}, {\em Embedded boundary algorithms and
  software for partial differential equations}, J. Phys., 125 (2008),
pp.~012084.

\bibitem{Karma98}
{\sc A.~Karma and W.-J. Rappel}, {\em Quantitative phase-field modeling of
  dendritic growth in two and three dimensions}, Physical Review E, 57 (1998),
pp.~4323--4349.

\bibitem{Kockelkoren03}
{\sc J.~Kockelkoren, H.~Levine, and W.~J. Rappel}, {\em Computational approach
  for modeling intra- and extracellular dynamics}, Phys. Rev., E 68 (2003),
p.~037702.

\bibitem{LL94}
{\sc R.J. Leveque and Z. Li}, {\em The immersed interface method for elliptic
  equations with discontinuous coefficients and singular sources}, SIAM J.
Numer. Anal., 31 (1994), pp.~1019-1044.

\bibitem{Levine05}
{\sc H.~Levine and W.~J. Rappel}, {\em Membrane-bound turing patterns},
Physical Review E, 72 (2005).

\bibitem{Li09}
{\sc X.~Li, J.~Lowengrub, A.~R\"atz, and A.~Voigt}, {\em Solving pdes in
  complex geometries: A diffuse-domain   approach}, Communications in
Mathematical Sciences, 7 (2009), pp.~81--107.

\bibitem{LiIto_2006}
{\sc Z. Li and K. Ito}, {\em The immersed interface method: Numerical solutions
  of {PDE}s involving interfaces and irregular domains}, SIAM Front. Appl.
Math., 33 (2006).

\bibitem{li12}
{\sc Z. Li, P. Song}, {\em An adaptive mesh refinement strategy for immersed
  boundary/interface methods}, Comm. Comput. Phys., 12 (2012), pp.~515-527.

\bibitem{lohner07}
{\sc R. Lohner and J.R. Cebral and F.F. Camelli and J.D. Baum and E.L. Mestreau
  and O.A. Soto}, {\em Adaptive embedded/immersed unstructured grid
  techniques}, Arch. Comput. Meth. Eng., 14 (2007), pp.~279-301.

\bibitem{lui09}
{\sc S.H. Lui}, {\em Spectral domain embedding for elliptic {PDE}s in complex
  domains}, J. Comput. Appl. Math., 225 (2009), pp.~541-557.

\bibitem{MacklinLowengrub_2006}
{\sc P. Macklin and J. Lowengrub}, {\em Evolving interfaces via gradients of
  geometry-dependent interior Poisson problems: Application to tumor growth},
J. Comput. Phys., 203 (2005), pp.~191-220.

\bibitem{MacklinLowengrub_2008}
{\sc P. Macklin and J. Lowengrub}, {\em A new ghost cell/level set method for
  moving boundary problems: Application to tumor growth}, J. Sci. Comput., 35
(2008), pp.~266-299.

\bibitem{OSK09}
{\sc M. Oevermann, C. Scharfenberg, and R. Klein}, {\em A sharp interface
  finite volume method for elliptic equations on Cartesian grids}, J. Comput.
Phys., 228 (2009), pp.~5184-5206.

\bibitem{osher03a}
{\sc S. Osher and R. Fedkiw}, {\em Level set methods and dynamic implicit
  surfaces}, Springer (2003).

\bibitem{osher88}
{\sc S. Osher and J.A. Sethian}, {\em Fronts propagating with
  curvature-dependent speed: Algorithms based on Hamilton-Jacobi formulations},
J. Comput. Phys., 79 (1988), pp.~12-49.

\bibitem{Papac10}
{\sc J. Papac, F. Gibou, and C. Ratsch}, {\em Efficient symmetric
  discretization for the Poisson, heat and Stefan-type problems with Robin
  boundary conditions}, J. Comput. Phys., 229 (2010), pp.~875-889.

\bibitem{Papac13}
{\sc J. Papac, A. Helgadottir, C. Ratsch, and F. Gibou}, {\em A level set
  approach for diffusion and Stefan-type problems with Robin boundary
  conditions on quadtree/octree adaptive Cartesian grids}, J. Comput. Phys.,
233 (2013), pp.~241-261.

\bibitem{Pego88}
{\sc R.~L. Pego}, {\em Front migration in the nonlinear {C}ahn-{H}illiard
  equation}, Proceedings of the Royal Society A, 422 (1988), pp.~261--278.

\bibitem{PRE11}
{\sc T. Preusser, M. Rumpf, S. Sauter, and L.O. Schwen}, {\em 3D composite
  finite elements for elliptic boundary value problems with discontinuous
  coefficients}, SIAM J. Sci. Comput., 35 (2011), pp.~2115-2143.

\bibitem{RamiereAngotBelliard_CMAME_2007}
{\sc I. Ramiere, P. Angot, and M. Belliard}, {\em A general fictitious domain
  method with immersed jumps and multilevel nested structured meshes}, J.
Comput. Phys., 225 (2007), pp.~1347-1387.

\bibitem{ratz}
{\sc A. R\"atz and A. Voigt}, {\em {PDE}s on surfaces---a diffuse interface
  approach}, Commun. Math. Sci., 4 (2006), pp.~575-590.

\bibitem{Reuter12}
{\sc M.~G. Reuter, J.~C. Hill, and R.~J. Harrison}, {\em Solving pdes in
  irregular geometries with multiresolution methods i: Embedded Dirichlet
  boundary conditions}, Computer Physics Communications, 183 (2012), pp.~1--7.

\bibitem{Sethian99}
{\sc J.A. Sethian}, {\em Level set methods and fast marching methods},
Cambridge University Press (1999), ISBN 0-521-64557-3.

\bibitem{SethianShan_2008}
{\sc J.A. Sethian and Y. Shan}, {\em Solving partial differential equations on
  irregular domains with moving interfaces, with applications to superconformal
  electrodeposition in semiconductor manufacturing}, J. Comput. Phys., 227
(2008), pp.~6411-6447.

\bibitem{Teigen09-b}
{\sc K.~E. Teigen, X.~Li, J.~Lowengrub, F.~Wang, and A.~Voigt}, {\em
  A diffuse-interface approach for modeling transport, diffusion and
  adsorption/desorption of material quantities on a deformable interface},
Communications in Mathematical Sciences, 7 (2009), pp.~1009--1037.

\bibitem{Teigen11}
{\sc K.E. Teigen, P. Song, A. Voigt, and J. Lowengrub}, {\em
  A diffuse-interface method for two-phase flows with soluble surfactants}, J.
Comput. Phys., 230 (2011), pp.~375-393.

\bibitem{Theillard13}
{\sc M. Theillard, L. F. Djodom, J.-L. Vié, and F. Gibou}, {\em A second-order
  sharp numerical method for solving the linear elasticity equations on
  irregular domains and adaptive grids – Application to shape optimization},
J. Comput. Phys., 233 (2013), pp.~430-448.

\bibitem{Uzgoren09}
{\sc E. Uzgoren, J. Sim, and W. Shyy}, {\em Marker-based, 3-D adaptive
  {C}artesian grid method for multiphase flows around irregular}, Comm.
Comput. Phys., 5 (2009), pp.~1-41.

\bibitem{wan12}
{\sc X.H. Wan and Z. Li}, {\em Some new finite difference methods for Helmholtz
  equations on irregular domains or with interfaces}, Disc. Cont. Dyn. Sys. B,
17 (2012), pp.~1155-1175.

\bibitem{Wise07}
{\sc S.~Wise, J.~Kim, and J.~Lowengrub}, {\em Solving the regularized, strongly
  anisotropic {C}ahn-{H}illiard equation by an adaptive nonlinear multigrid
  method}, Journal of Computational Physics, 226 (2007), pp.~414--446.

\bibitem{xia11}
{\sc K. Xia, M. Zhan, and G. Wei}, {\em MIB method for elliptic equations with
  multimaterial interfaces}, J. Comput. Phys., 230 (2011), pp.~4588-4615.

\bibitem{zhao09}
{\sc S. Zhao and G. Wei}, {\em Matched interface and boundary (MIB) for the
  implementation of boundary conditions in high-order central finite
  differences}, Int. J. Numer. Meth. Eng., 77 (2009), pp.~1690-1730.

\bibitem{zhao10}
{\sc S. Zhao}, {\em High order matched interface and boundary methods for the
  Helmholtz equation in media with arbitrarily curved interfaces}, J. Comput.
Phys., 229 (2010), pp.~3155-3170.

\bibitem{zhong07}
{\sc X.L. Zhong}, {\em A new high-order immersed interface method for solving
  elliptic equations with imbedded interface of discontinuity}, J. Comput.
Phys., 225 (2007), pp.~1066-1099.

\bibitem{zhou06}
{\sc Y.C. Zhou, S. Zhao, M. Feig, and G.W. Wei}, {\em High order matched
  interface and boundary method for elliptic equations with discontinuous
  coefficients and singular sources}, J. Comput. Phys., 213 (2006), pp.~1-30.

\bibitem{zhou12}
{\sc Y.C. Zhou, J. Liu, and D.L. Harry}, {\em A matched interface and boundary
  method for solving multiflow Navier-Stokes equations with applications to
  geodynamics}, J. Comput. Phys., 231 (2012), pp.~223-242.

\bibitem{Zhu12}
{\sc Y. Zhu, Y. Wang, J. Hellrung, A. Cantarero, E. Sifakis, and J. M. Teran},
{\em A second-order virtual node algorithm for nearly incompressible linear
  elasticity in irregular domains}, J. Comput. Phys., 231 (2012),
pp.~7092-7117.

\end{thebibliography}

\end{document}